\documentclass[10pt]{article}
\usepackage{amsmath,amsfonts,amscd,amssymb,hyperref}
\DeclareMathAlphabet\EuR{U}{eur}{m}{n}
\SetMathAlphabet\EuR{bold}{U}{eur}{b}{n}

\begin{document}

\newcommand{\bigtit}{\begin{bf} \begin{center} \begin{LARGE}
``The p-chain spectral sequence''
\\[3mm] by\\[3mm] James F. Davis and Wolfgang L\"uck\\[4mm] 
\end{LARGE}\end{center}\end{bf}
\bigskip}

\newcommand{\cala}{{\mathcal A}}
\newcommand{\calb}{{\mathcal B}}
\newcommand{\calc}{{\mathcal C}}
\newcommand{\cald}{{\mathcal D}}
\newcommand{\cale}{{\mathcal E}}
\newcommand{\calf}{{\mathcal F}}
\newcommand{\calg}{{\mathcal G}}
\newcommand{\calh}{{\mathcal H}}
\newcommand{\calm}{{\mathcal M}}
\newcommand{\caln}{{\mathcal N}}
\newcommand{\calp}{{\mathcal P}}
\newcommand{\caltr}{{\mathcal T}\!{\mathcal R}}
\newcommand{\cals}{{\mathcal S}}
\newcommand{\calz}{{\mathcal Z}}

\newcommand{\calmfin}{{\mathcal M}\!{\mathcal F}\!{\mathcal I}\!{\mathcal N}}
\newcommand{\calfin}{{\mathcal F}\!{\mathcal I}\!{\mathcal N}}
\newcommand{\calvc}{{\mathcal V}\!{\mathcal C}}
\newcommand{\calall}{{\mathcal A}\!{\mathcal L}\!{\mathcal L}}
\newcommand{\calsub}{{\mathcal S}{\mathcal U}\!{\mathcal B}}
\newcommand{\inc}{\operatorname{inc}}

\newcommand{\da}{\downarrow}

\newcommand{\rr}{{\mathbb R}}
\newcommand{\hh}{{\mathbb H}}
\newcommand{\zz}{{\mathbb Z}}
\newcommand{\nn}{{\mathbb N}}
\newcommand{\zmod}{\zz/2}
\newcommand{\cc}{{\mathbb C}}
\newcommand{\qq}{{\mathbb Q}}
\newcommand{\ff}{{\mathbb F}}

\newcommand{\bfA}{\ensuremath{\mathbf{A}}}
\newcommand{\bfB}{\ensuremath{\mathbf{B}}}
\newcommand{\bfE}{\ensuremath{\mathbf{E}}}
\newcommand{\bfF}{\ensuremath{\mathbf{F}}}
\newcommand{\bfK}{\ensuremath{\mathbf{K}}}
\newcommand{\bfL}{\ensuremath{\mathbf{L}}}
\newcommand{\bfT}{\ensuremath{\mathbf{T}}}
\newcommand{\bfS}{\ensuremath{\mathbf{S}}}
\newcommand{\bfDelta}{\ensuremath{\mathbf{\Delta}}}
\newcommand{\boldDelta}{\mbox{\boldmath $\Delta$}}
\newcommand{\boldDeltasmall}{\mbox{{\scriptsize \boldmath $\Delta$}}}

\newcommand{\alg}{\operatorname{alg}}
\newcommand{\aut}{\operatorname{aut}}
\newcommand{\cok}{\operatorname{cok}}
\newcommand{\colim}{\operatorname*{colim}}
\newcommand{\cone}{\operatorname{cone}}
\newcommand{\cyl}{\operatorname{cyl}}
\newcommand{\Ext}{\operatorname{Ext}}
\newcommand{\func}{\operatorname{functor}}
\newcommand{\holim}{\operatorname*{holim}}
\newcommand{\hocolim}{\operatorname*{hocolim}}
\newcommand{\hoinvlim}{\operatorname*{hoinvlim}}
\newcommand{\id}{\operatorname{id}}
\newcommand{\im}{\operatorname{im}}
\newcommand{\Is}{\operatorname{Is}}
\newcommand{\invlim}{\operatorname*{invlim}}
\newcommand{\map}{\operatorname{map}}
\newcommand{\mor}{\operatorname{mor}}
\newcommand{\Ob}{\operatorname{Ob}}
\newcommand{\op}{\operatorname{op}}
\newcommand{\pr}{\operatorname{pr}}
\newcommand{\topo}{\operatorname{top}}
\newcommand{\Tor}{\operatorname{Tor}}
\newcommand{\trans}{\operatorname{trans}}
\newcommand{\Wh}{\operatorname{Wh}}

\let\sect=\S
\newcommand{\curs}{\EuR}
\newcommand{\CATEGORIES}{\curs{CATEGORIES}}
\newcommand{\CHAIN}{\curs{CHAIN}}
\newcommand{\COMPLEXES}{\curs{COMPLEXES}}
\newcommand{\MOD}{\curs{MOD}}
\newcommand{\Or}{\operatorname{Or}}
\newcommand{\SETS}{\curs{SETS}}
\newcommand{\SPACES}{\curs{SPACES}}
\newcommand{\SPECTRA}{\curs{SPECTRA}}


\newtheorem{theorem}{Theorem}[section]
\newtheorem{proposition}[theorem]{Proposition}
\newtheorem{lemma}[theorem]{Lemma}
\newtheorem{definition}[theorem]{Definition}
\newtheorem{example}[theorem]{Example}
\newtheorem{diagram}[theorem]{Diagram}
\newtheorem{remark}[theorem]{Remark}
\newtheorem{notation}[theorem]{Notation}
\newtheorem{refnumber}[theorem]{}
\newtheorem{corollary}[theorem]{Corollary}
\newtheorem{assumption}[theorem]{Assumption}
\newtheorem{conjecture}[theorem]{Conjecture}

{\catcode`@=11\global\let\c@equation=\c@theorem}
\renewcommand{\theequation}{\thetheorem}

\renewcommand{\theenumi}{\alph{enumi}}
\renewcommand{\labelenumi}{(\theenumi)}

\newcommand{\tit}[2]{\begin{bf} \begin{center} \begin{Large}
\section{#1} \label{sec: #2}
\end{Large}\end{center}\end{bf}}


\newcommand{\comsquare}[8]{
\begin{center}
$\begin{CD} #1 @>#2>> #3\\ @V{#4}VV @VV{#5}V\\ #6 @>>#7> #8
\end{CD}$
\end{center}}

\newcommand{\squarematrix}[4]{\left( \begin{array}{cc} #1 & #2 \\ #3 &
#4
\end{array} \right)}


\newcommand{\cata}[3]{\, #1 \!\downarrow\! #2 \!\downarrow #3\,}
\newcommand{\catb}[2]{\, #1 \!\downarrow\! #2 \,}


\newcommand{\point}{\star}

\newcommand{\proofl}{{\bf \underline{\underline{Proof}}: }}
\newcommand{\qedl}{\hspace{10mm} \rule{2mm}{2mm}\bigskip}


\bigtit

\typeout{-----------------------  Abstract  ------------------------}
\begin{abstract}
We introduce a new spectral sequence called the $p$-chain spectral
sequence which converges to the (co-)homology of a contravariant
$\calc$-space with coefficients in a covariant $\calc$-spectrum
for a small category $\calc$. It is different from the corresponding
Atiyah-Hirzebruch type spectral sequence. It can be used in
combination with the Isomorphism Conjectures of Baum-Connes and
Farrell-Jones to compute algebraic $K$- and $L$-groups of group rings
and topological $K$-groups of reduced group $C^*$-algebras.
\smallskip


\end{abstract}


\typeout{-----------------------  Introduction ------------------------}

\setcounter{section}{-1}
\tit{Introduction}{Introduction}

In \cite{Davis-Lueck (1998)} we defined abelian groups $H^\calc_n(X;\bfE)$ 
when $X \colon  \calc \to \SPACES$ is a contravariant
functor and $\bfE \colon  \calc \to \SPECTRA$ is a covariant functor
\footnote{We also defined cohomology groups $H^n_\calc(X;\bfE)$ when
$X \colon \calc \to \SPACES$ and $\bfE \colon  \calc \to \SPECTRA$ are 
{\em both} contravariant functors.  In this introduction we will only
discuss homology in order to simplify the exposition.}.  
These abelian groups behave like a generalized homology
theory defined on $\calc$-spaces, for example, there is a long exact
Mayer-Vietoris  sequence.  They are weak homotopy
invariant; given a map $f \colon  X \to Y$ of $\calc$-spaces 
(i.e. a natural transformation) which induces a weak homotopy
equivalence $f(c) \colon  X(c) \to Y(c)$ for all objects $c \in \Ob \calc$, 
then $f_* \colon  H_n^\calc(X;\bfE) \to H^\calc_n(Y;\bfE)$
is an isomorphism.  

Three special cases will illustrate these groups:
\begin{enumerate}
\item  Fixing an object $c \in \Ob \calc$, 
$H^\calc_n(\mor_\calc(?,c)) = \pi_n(\bfE(c))$.  This should be 
thought of as giving the coefficients of the generalized
homology theory;

\item  If $X = \point $ is the constant functor, then 
$H_n^\calc(\point;\bfE) = \pi_n(\hocolim_\calc \bfE)$;

\item  If $\calc$ is a category with a single object, all of 
whose morphisms are isomorphisms, our generalized
 homology theory reduces to Borel homology.  More precisely, 
let $G$ be the group of morphisms and let $X$ be a 
$CW$-complex with an action of $G$ by cellular maps, then
$$
H_n^\calc(X;\bfE) = H_n^G(X;\bfE) =: \pi_n((X \times EG)_+ \wedge_G \bfE).
$$ 
\end{enumerate}

In \cite{Davis-Lueck (1998)} we gave a spectral sequence 
converging to $H_{p+q}^\calc(X;\bfE)$ whose $E^2$-term is
$E^2_{p,q} = H_p^\calc(X; \pi_q(\bfE))$.
This spectral sequence is both quite useful and quite 
standard.  It is related to spectral sequences of
Atiyah-Hirzebruch, Bousfield-Kan, and Quinn.  The point of 
the current paper is to introduce a {\em new} spectral sequence. 
It converges to $H_{p+q}^\calc(X;\bfE)$ and has as $E^1$-term
$$E^1_{p,q} ~ = ~ 
\coprod_{\substack{\overline{c}_{\ast} \in \prod_{i=0}^p \Is( \calc ),\\
S(\overline{c}_*) \not= \emptyset}}~
 H_q^{\calc}(X(c_p) \times_{\aut(c_p)}
S(\overline{c}_{\ast})
\times_{\aut(c_0)} \mor(?,c_0);\bfE),$$
where $\Is(\calc)$ is the set of isomorphism classes of 
objects in $\calc$ and $S(\overline{c}_{\ast})$ is a certain
$\aut(c_p)$-$\aut(c_0)$-set.  
In the special case of an orbit category
$\Or(G,\calf)$, where $G$ is a (discrete) group and $\calf$ is a family 
of subgroups satisfying the condition  $gHg^{-1}
\subseteq H \Rightarrow gHg^{-1}
= H$ for $g \in G$ and $H \in \calf$, then $E^1_{p,q}$ can be 
indexed by ``$p$-chains"
$$
(H_0) < (H_1) < \ldots <(H_p)
$$
where $(H_i)$ is the conjugacy class of $H_i \in \calf$ 
and $(H_i) < (H_{i+1})$ means that there exists $g \in G$ so that
$gH_ig^{-1} \subsetneqq H_{i+1}$.  More generally, if $\calc$ 
is a EI-category (endomorphisms are isomorphisms), then $E^1_{p,q}$
is indexed by similar $p$-chains.  The $E^1_{p,q}$-term looks 
formidable, but we discuss many simplifications and reinterpretations
in cases of interest.  In general, the $E^1$-term of the spectral 
sequence is related to Borel homology of the groups
$\{\aut(c)\}_{c
\in \Ob \calc}$ and the differentials are related to assembly maps.

It should be emphasized that the motivation for this 
abstract-looking spectral sequence comes from geometric topology (the
surgery classification of manifolds), geometry (which manifolds 
admits metrics of positive scalar curvature), and analysis (the
study of $K$-theory of group $C^*$-algebras), in conjunction 
with the study of (fundamental) groups which are infinite, but
contain torsion. The connection with these various subjects 
comes through isomorphism conjectures and assembly maps, for more on
this see
\cite{Davis-Lueck (1998)} and the end of Section 
\ref{sec: Review of Spaces over a Category and Assembly Maps} of the current paper.

In Section \ref{sec: Review of Spaces over a Category and Assembly Maps}
 we review the definition of $H^\calc_n(X;\bfE)$ 
and discuss maps connected with a functor $F \colon  \calc \to \cald$.  In
Section \ref{sec: The p-Chain Spectral Sequence}
 we derive the $p$-chain spectral sequence and its 
differentials and discuss simplifications which occur when the
category $\calc$ is left-free.  In Section 
\ref{sec: Assembly Maps} we discuss  the important special case
$H^\calc_n(\point;\bfE)$ where $\calc$ is the restricted orbit category.  Here, by using the 
Cofinality Theorem and an analogue of Quillen's Theorem A, computations can be simplified.
We also discuss differentials in the $p$-chain spectral sequence, which 
often turn out to be assembly maps themselves. In Section \ref{sec: Examples} we
give examples of groups where the methods of the preceeding sections in
combination with the Isomorphism Conjectures of Baum-Connes and
Farrell-Jones lead to explicit computations of 
algebraic $K$- and $L$-groups of group rings
and topological $K$-groups of reduced group $C^*$-algebras.

The $p$-chain spectral sequence is a generalization of the spectral 
sequence in \cite[Chapter 17]{Lueck (1989)} and is related to
the paper of Slominska \cite{Slominska (1989)}.

The authors thank Bill Dwyer and Reiner Vogt for useful conversations.  The second author
thanks the Max-Planck Institut f\"ur Mathematik in Bonn
for its hospitality
during his stay in November/December 2002 when parts of this paper were written.

The paper is organized as follows:\\[1mm]
\begin{tabular}{ll}  \ref{sec: Introduction}. & Introduction
\\  
\ref{sec: Review of Spaces over a Category and Assembly Maps}. & Review of Spaces over a
Category and Assembly maps
\\ 
\ref{sec: The p-Chain Spectral Sequence}. & The $p$-Chain Spectral Sequence
\\ 
\ref{sec: Assembly Maps}. & Assembly Maps
\\ 
\ref{sec: Examples}. & Examples
\\  
 & References
\end{tabular}

\bigskip


\typeout{-------------------------- section 1 --------------------------}

\tit{Review of Spaces over a Category and Assembly Maps}
{Review of Spaces over a Category and Assembly Maps}

In this section we review some basic facts from
\cite{Davis-Lueck (1998)} for the convenience of the reader.

Let $\SPACES$ and  $\SPACES_+$ be the categories of topological spaces
and pointed topological spaces respectively. We will
always work in the category of compactly generated spaces (see
\cite{Steenrod (1967)} and
\cite[I.4]{Whitehead (1978)}). A {\em spectrum}
$$\mathbf{E} = \{(E(n),\sigma(n)) \mid n \in \zz\}$$
is a sequence of pointed spaces
$\{E(n) \mid n \in \zz\}$ together with pointed maps
$\sigma(n) \colon  E(n) \wedge S^1 \to E(n+1)$, called {\em
structure maps}. A {\em (strong) map} of
spectra
$\mathbf{f} \colon  \mathbf{E} \to \mathbf{E}^{\prime}$ is a sequence of maps
$f (n) \colon  E(n) \to E^{\prime}(n)$ which are compatible
with the structure maps, i.e.
$f(n+1) \circ \sigma(n) ~ = ~
\sigma^{\prime}(n) \circ \left(f(n) \wedge \id_{S^1}\right)$ holds for all
$n \in \zz$. This should not be
confused with the notion of map of spectra in the stable category (see
\cite[III.2.]{Adams (1974)}). A spectrum is
called an {\em $\Omega$-spectrum} if the adjoint of each structure map
$E(n) \to \Omega E(n+1)$ is a weak homotopy
equivalence. The homotopy groups of a spectrum are
defined by
$$\pi_i(\mathbf{E}) ~ = ~
\colim_{k \to \infty} \pi_{i+n}(E(n))$$ where the system $\pi_{i+n}(E(n))$
is given by the composite
$$\pi_{i+n}(E(n)) \xrightarrow{S} \pi_{i+n+1}(E(n)\wedge S^1)
\xrightarrow{\sigma (n)_*} \pi_{i+n+1}(E(n +1))$$ of the suspension
homomorphism and the homomorphism induced by
the structure map. 

Let $\calc$ be a small category, i.e. a category such that the objects and
the morphisms form sets. A {\em
covariant (or contravariant) $\calc$-space, pointed
$\calc$-space,  $\calc$-spectrum, $\ldots$ } is a covariant (or
contravariant) functor from $\calc$ to
$\SPACES$, $\SPACES_+$, $\SPECTRA$, $\ldots $ and a morphism is a
natural transformation. 

Let $X$ be a contravariant and $Y$ be a covariant $\calc$-space. Define
their {\em tensor product} to be the space
$$X \otimes_{\calc} Y ~ = ~
\coprod_{c \in \Ob(\calc)} X(c) \times Y(c)/\sim$$ where $\sim$ is the
equivalence relation which is generated by
$(x\phi,y) \sim (x,\phi y)$ for all morphisms
$\phi  \colon  c \to d$ in $\calc$ and points $x \in X(d)$ and
$y \in Y(c)$. Here $x\phi$ stands for
$X(\phi)(x)$ and $\phi y$ for $Y(\phi )(y)$. If $X$ and $Y$ are $\calc$-spaces of the same
variance,  denote by
$\hom_{\calc}(X,Y)$ the space of maps of $\calc$-spaces from $X$
to $Y$ with the subspace topology coming
from the obvious inclusion into $\prod_{c \in \Ob(\calc)} \map(X(c),Y(c))$. 
If $Y$ is a $\calc$-space and $Z$ a space, let $Y \times Z$ be the
$\calc$-space with the  same
variance as $Y$ whose value at the object $c$ in $\calc$ is
 $Y(c) \times Z$. For $\calc$-spaces $X$ and $Y$ of the same variance, the set of
homotopy classes of maps of $\calc$-spaces $[X,Y]^\calc$ is defined using maps of
$\calc$-spaces $X
\times [0,1]
\to Y$.  If
$Y$ is a
$\calc$-space and
$Z$ a space, let $\map(Y,Z)$  be the
$\calc$-space with the opposite 
variance as $Y$ whose value at the object $c$ in $\calc$ is
$\map(Y(c),Z)$. If $X$ is a
contravariant $\calc$-space, $Y$ a covariant
$\calc$-space and
$Z$ a space, there is a canonical adjunction homeomorphism
\cite[Lemma 1.5]{Davis-Lueck (1998)}
$$\map(X \otimes_{\calc} Y,Z) \xrightarrow{\cong}
\hom_{\calc}(X,\map(Y,Z)).$$ 

All of the above notions also make sense for pointed
spaces; one has to substitute wedges for
disjoint unions, smash products for cartesian products, and pointed
mapping spaces for mapping spaces.

A {\em contravariant $\calc$-$CW$-complex}
$X$ is a contravariant
$\calc$-space $X$ together with a filtration
$$\emptyset = X_{-1} \subseteq X_0 \subseteq X_1 \subseteq X_2 \subseteq \ldots
\subseteq X_n \subseteq \ldots \subseteq X = \bigcup_{n \ge 0} X_n$$ such that
$X = \colim_{n \to \infty} X_n$ and for
any $n \ge 0$ the {\em $n$-skeleton} $X_n$ is obtained from the
$(n-1)$-skeleton $X_{n-1}$ by attaching $\calc$-$n$-cells, i.e. there
exists a pushout of $\calc$-spaces of
the form
$$\begin{CD}
\coprod_{i \in I_n} \mor_{\calc}(?,c_i) \times S^{n-1} @>>> X_{n-1} \\
@VVV @VVV \\ {\coprod_{i \in I_n}
\mor_{\calc}(?,c_i) \times D^n} @>>> {X_n}
\end{CD}$$ where the vertical maps are inclusions, $I_n$ is an index set,
and the $c_i$ are objects of $\calc$. 
The definition of a {\em covariant
$\calc$-$CW$-complex} is analogous. In \cite{Davis-Lueck (1998)} these were called
contravariant free and covariant free $\calc$-$CW$-complexes, we will omit the word free here.
One of the main properties of $CW$-complexes carries
over to $\calc$-$CW$-complexes, namely, a map $f\colon Y \to Z$
of $\calc$-spaces is a {\em weak homotopy
equivalence}, i.e. $f(c)$ is a weak homotopy equivalence of spaces for all
objects $c$ in $\calc$, if and only if
for any $\calc$-$CW$-complex $X$ the induced map
$$f_{\ast}\colon  [X,Y]^{\calc} \to [X,Z]^{\calc} \hspace{10mm} [g]
\mapsto [g \circ f]$$ between the
homotopy classes of maps of $\calc$-spaces is bijective \cite[Theorem 3.4]{Davis-Lueck (1998)}.
In particular Whitehead's Theorem carries over: a
map of  $\calc$-$CW$-complexes is a homotopy equivalence if and
only if it is a weak homotopy equivalence. A
$\calc$-$CW$-approximation $(X,f)$ of a $\calc$-space $Y$
consists of a  $\calc$-$CW$-complex $X$
together with a weak homotopy equivalence $f\colon  X \to Y$. Such a
$\calc$-$CW$-approximation always exists, there is even a functorial
construction. If $(X,f)$ and 
$(X^{\prime},f^{\prime})$ are two $\calc$-$CW$-approximations of a $\calc$-space
$Y$,
there is a homotopy equivalence $h\colon  X \to X^{\prime}$ of
$\calc$-spaces which is determined uniquely up to homotopy by the property that
$f^{\prime} \circ h$ and $f$ are homotopic. 

Let $(X,A)$ be a pair of contravariant pointed $\calc$-spaces. Denote the reduced cone of
the pointed space $A$ by $\cone (A)$.
For a covariant $\calc$-spectrum $\mathbf{E}$  define
$$\mathbf{E}^{\calc}_q(X,A)  ~ = ~
\pi_q(X \cup_{A} \cone (A) \otimes_{\calc} \mathbf{E}).$$ Given a  contravariant
$\calc$-spectrum $\mathbf{E}$, define
$$\mathbf{E}^q_{\calc}(X,A) ~ = ~
\pi_{-q}(\hom_{\calc}(X \cup_A \cone (A),\mathbf{E})).$$ 
Define $(X_{+},A_{+})$ to be the pair of contravariant pointed $\calc$-spaces which is
obtained from $(X,A)$ by adding a
disjoint base point. Let $(u,v) \colon  (X^{\prime},A^{\prime}) ~
\to (X,A)$ be a
$\calc$-$CW$-approximation.  For a covariant $\calc$-spectrum $\mathbf{E}$  define the
{\em homology of $(X,A)$ with coefficients in $\mathbf{E}$} by
$$H^{\calc}_q(X,A;\mathbf{E})  ~ = ~
\mathbf{E}^{\calc}_q(X^{\prime}_{+},A^{\prime}_{+}).$$

Given a contravariant $\calc$-spectrum $\mathbf{E}$,
define the {\em cohomology of $(X,A)$ with
coefficients in $\mathbf{E}$} by
$$H_{\calc}^q(X,A;\mathbf{E})  ~ = ~
\mathbf{E}_{\calc}^q(X^{\prime}_{+},A^{\prime}_{+}).$$ When $A$ is empty
we omit it from the notation.

Then $H_q^{\calc}(X,A;\mathbf{E})$ and
$H^q_{\calc}(X,A;\mathbf{E})$ are unreduced homology and cohomology
theories on pairs of $\calc$-spaces which
satisfy the  WHE-axiom, which says that a weak homotopy
equivalence induces an isomorphism on (co-)homology. The homology theory satisfies the
disjoint union axiom. The cohomology theory satisfies
the disjoint union axiom provided that $\bfE$ is a $\calc$-$\Omega$-spectrum.  Let $\point_\calc$,
or briefly $\point$, be the $\calc$-space which takes each object to a
point.  Then for $(X,A) = (\point, \emptyset)$ then the above notions
reduce to 
$$\pi_q(\hocolim_{\calc} \bfE)
$$ and 
$$\pi_{-q}(\holim_{\calc} \bfE)
$$ respectively.  If $\calc$ is the category associated to a group $G$, i.e. $\calc$ has  a single
object, the morphisms in $\calc$ are in one-to-one correspondence with $G$, and the  composition law 
in $\calc$ corresponds to multiplication in $G$,
then $H_q^{\calc}(X,A;\mathbf{E})$  can be identified with Borel homology  $\pi_q(EG \times
(X_+,A_+)
\wedge_G \bfE)$, and similarly for cohomology, provided that
$(X,A)$ is a $CW$-pair with a $G$-action by cellular maps.  Finally, if
$\calc$ is the trivial category with precisely one object and morphism,
these notions reduce to the standard
notions of the homology and cohomology of spaces given by a
spectrum.

Given a functor $F \colon  \calc \to \cald$  and a $\cald$-space (or spectrum) $Y$, define the {\em
restriction of 
$Y$ with respect to $F$} to be the $\calc$-space (or spectrum) $F^*Y(c) = Y(F(c))$.  For $X$ a
contravariant, respectively covariant,
$\calc$-space define the {\em induction of $X$ with respect to $F$} to be the $\cald$-space
$$
F_*X(??) = X(?) \otimes_\calc \mor_\cald(??,F(?)),
$$
respectively
$$
F_*X(??) = \mor_\cald(F(?),??) \otimes_\calc X(?).
$$
There are natural adjunction isomorphisms (see \cite[Lemma 1.9]{Davis-Lueck (1998)})
\begin{eqnarray}
X\otimes_\calc F^*Y & \cong  & F_*X \otimes_\cald Y, \label{adjunction (F_*,F^*) for otimes}
\\
\hom_\calc (X,F^*Y) & \cong  & \hom_\cald(F_*X,Y).\label{adjunction (F_*,F^*) for hom}
\end{eqnarray}

\begin{lemma}\label{lem: induction}
Let $F \colon  \calc \to \cald$ be a covariant functor.  
\begin{enumerate}
\item \label{lem: induction: cofibrations and fibrations}

If $f \colon  X \to Y$ is a fibration of $\cald$-spaces, then  
$F^*f \colon  F^*X \to F^*Y$ is a fibration of $\calc$-spaces.  If $f \colon  X
\to Y$ is a cofibration of
$\calc$-spaces, then $F_*f\colon  F_*X \to F_*Y$ is a cofibration of
$\cald$-spaces;

\item\label{lem: induction: CW-complexes}
If $X$ is a  $\calc$-$CW$-complex, then $F_*X$ is a  $\cald$-$CW$-complex;

\item\label{lem: induction: map Phi_F}
Let $X$ be a contravariant $\calc$-space and $\bfE$ be a covariant $\cald$-spectrum.  
Then there is a map, natural in $X$
and $\bfE$,
$$
\Phi_ F:H_q^\calc(X; F^*\bfE) \to H_q^\cald(F_*X;\bfE).
$$
If either 
\begin{enumerate}
\item[(a)]
$X$ is a  $\calc$-$CW$-complex,  or,

\item[(b)] 
for all objects $d$ of $\cald$ the covariant $\calc$-set $\mor_\cald(d,F(?))$  is isomorphic
to a disjoint union of covariant $\calc$-sets of the form $\mor_\calc(c_0,?))$,
\end{enumerate}
then $\Phi_ F:H_q^\calc(X; F^*\bfE) \to H_q^\cald(F_*X;\bfE)$ is 
an isomorphism. 

Similar statements hold true for cohomology.
 \end{enumerate}

\end{lemma}

\proofl 
\eqref{lem: induction: cofibrations and fibrations}  
To show that $F^*f$ is a fibration, one sets up the 
homotopy lifting problem for $\calc$-spaces, solves the adjoint
problem for
$\cald$-spaces, and uses the adjoint property to translate the 
solution back to $\calc$-spaces.  The proof for
cofibrations is similar.
\\[2mm]
\eqref{lem: induction: CW-complexes}
This follows from two facts.  First, for any $c \in \Ob C$, one can  identify
$F_*\mor_\calc(?,c) \cong \mor_\cald(?,F(c))$.  Second, since $F_*$ 
has a right adjoint, it commutes with colimits, in particular, it commutes with pushouts.
\\[2mm]
\eqref{lem: induction: map Phi_F}
Define the map $\Phi_F$ so that when $X' \to X$ is a $\calc$-$CW$-approximation, 
the following diagram commutes
$$
\begin{CD}
(F^*\bfE)^\calc_q(X') @>\cong>> (\bfE)^\cald_q(F_* X') \\
@V\cong VV 
@VV\cong V \\
H^\calc_q(X';  F^*\bfE)  @. H^\cald_q(F_*X';  \bfE)\\
@V\cong VV @VVV \\
H^\calc_q(X;  F^*\bfE) @>\Phi_F>> H^\cald_q(F_*X;  \bfE).
\end{CD}
$$
Uniqueness of $\calc$-$CW$-approximations (up to homotopy) gives that $\Phi_F$ is well-defined.  
If condition (a) holds, then the identity map 
$X \to X$ is a $\calc$-$CW$-approximation and the claim follows.

If condition (b) is satisfied, then a weak homotopy equivalence
 $X' \to X$ gives a weak homotopy equivalence
$ F_*X' \to F_*X $ and hence the lower right vertical
map in the above diagram is an isomorphism.
\qedl

For the applications we are most interested in the most important example is the {\em orbit
category} $\Or(G)$ defined for a group $G$.  The objects are homogeneous $G$-spaces $G/H$ and the
morphisms are the $G$-maps.  More generally, for a family $\calf$ of subgroups of $G$, define
the \emph{restricted orbit category} $\Or(G,\calf)$ 
to be the category whose objects are the homogeneous
$G$-spaces $G/H$ where $H \in \calf$ and the morphisms are the $G$-maps.  Some examples for
$\calf$ are $\caltr$,  $\calfin$, $\calvc$, and ${\cal ALL}$, 
which are the families consists of the
trivial group, the finite subgroups, the virtually cyclic subgroups, and all subgroups
respectively.

The remainder of this section is not necessary for the 
discussion of the $p$-chain spectral sequence itself in Section 
\ref{sec: The p-Chain Spectral Sequence}, but does give the motivation 
for this paper and is necessary for Section \ref{sec: Assembly Maps}.  We will review the point of
view of \cite{Davis-Lueck (1998)} concerning assembly maps and 
the Farrell-Jones and Baum-Connes Isomorphism Conjectures.  

\begin{definition}  \label{def: assembly map associated to a functor}
Let $F \colon  \calb \to \calc$ and $\bfE \colon  \calc \to \SPECTRA$ be covariant functors. 
Then the composite
$$
H^\calb_q(\point;F^*\bfE) \xrightarrow{\Phi_F}  H^\calc_q(F_*\point;\bfE) 
\xrightarrow{ H^\calc_q(\pr;\bfE)}  H^\calc_q(\point;\bfE)
$$
is called the {\em assembly map induced by $F$}, where $\Phi_F$ is the map appearing in
Lemma  \ref{lem: induction} \eqref{lem: induction: map Phi_F} 
and $\pr\colon F_*\point \to \point$ is the constant map at each object.
We sometimes abbreviate  $F^*\bfE$ by $\bfE$.
\end{definition}

\begin{remark} 
\em This map can be identified with the map 
$\pi_q\left(\hocolim_\calb F^*\bfE\right) \to \pi_q\left(\hocolim_\calc \bfE\right)$ 
induced by the functor $F$.\em
\end{remark}

The three $\Or(G)$-spectra which are useful to us are covariant functors 
\begin{align*}
&\bfK^{\alg} \colon  \Or(G) \to \SPECTRA \\
&\bfL^{\langle j \rangle} \colon  \Or(G) \to \SPECTRA \\
&\bfK^{\topo} \colon  \Or(G) \to \SPECTRA. \\
\end{align*} 
These functors were constructed in \cite[Section 2]{Davis-Lueck (1998)}, but there was a
problem with the construction of $\bfK^{\topo}$ connected with the pairing  on
\cite[p. 217]{Davis-Lueck (1998)}.  
This problem can easily be fixed and the construction can be replaced by  more refined
ones (see \cite{Joachim (2001)}).  The key property of these
functors is that $\pi_q(\bfK^{\alg}(G/H)) = K_q(RH)$ for a fixed ring $R$,
 $\pi_q(\bfL^{\langle -\infty \rangle}(G/H)) = L_q(RH)$ for a fixed 
ring $R$ with involution, and 
$\pi_q(\bfK^{\topo}(G/H)) = K_q(C^*_rH)$. Here  $K_q(C^*_rH)$ is the
$K$-theory of the real or complex reduced $C^*$-algebra of $H$.  
In connection with $L$-groups we use the involution on $RG$
sending $r \cdot g$ to $\overline{r} \cdot g^{-1}$
The index $j \in \zz \cup {-\infty}$ on the $L$-theory is
the $K$-theory decoration; the important cases for us are $j = -\infty$ which
arises in the isomorphism conjecture, and the case $j=2$ which is used in Wall's
book \cite{Wall (1970)} to classify manifolds. The Isomorphism Conjecture of Farrell and Jones is
not true for the decorations $ j = 0,1,2$, which correspond to  the decorations $p$, $h$ and $s$
appearing in the literature \cite{Farrell-Jones-Lueck (2002)}. 

The Isomorphism Conjecture of Baum-Connes  for a group $G$ is that the assembly map 
associated to the inclusion functor $I \colon \Or(G,\calfin) \to \Or(G)$ 
(see Definition \ref{def: assembly map associated to a functor}) yields
an isomorphism
$$H_q^{\Or(G,\calfin)}(\point;\bfK^{\topo}) \to 
H_q^{\Or(G)}(\point;\bfK^{\topo}) = K_q(C^*_r(G)).
$$
The Isomorphism Conjecture of Farrell-Jones for $RG$ says that 
$$
H_q^{\Or(G,\calvc)}(\point;\bfK^{\alg}) \to 
H_q^{\Or(G)}(\point;\bfK^{\alg}) = K_q(RG)
$$
and 
$$
H_q^{\Or(G,\calvc)}(\point;\bfL^{\langle -\infty\rangle}) \to
H_q^{\Or(G)}(\point;\bfL^{\langle -\infty\rangle}) = L^{\langle -\infty \rangle}_q(RG)
$$
are isomorphisms. We will mainly consider the case $R = \zz$.
The point of these conjectures is 
that they express the target, which is the group one 
wants to compute, by the source, which only involves the $K$-theory 
of the family of finite or virtually cyclic subgroups and is much easier to compute.
In the case of the family $\calfin$ the source can  
rationally be computed by equivariant Chern characters
\cite{Lueck (2002b)}, \cite{Lueck (2002c)}. The $p$-chain spectral sequence is an important tool
for integral computations which are much harder. More information about these conjectures
can be found for instance in 
\cite{Baum-Connes-Higson (1994)},
\cite{Farrell-Jones (1993a)}, 
\cite{Mislin (2002)}
and \cite{Valette (1997)}.

\bigskip


\typeout{-------------------------- section 2 --------------------------}

\tit{The $p$-Chain Spectral Sequence}{The p-Chain Spectral Sequence}

  We establish a spectral sequence converging to the homology (respectively
cohomology)
 of a  $\calc$-space $X$ with coefficients in a  $\calc$-spectrum
$\bfE$. In the special case where $X = \point$, this gives a spectral sequence converging
to $\pi_*(\hocolim_\calc \bfE)$ (respectively $\pi_*(\holim_\calc \bfE)$).  It is different
from the standard spectral sequence
\cite[Theorem 4.7]{Davis-Lueck (1998)},
\cite[Theorem 8.7]{Quinn (1982a)}, \cite[XII, 5.7 on page 339 and XI, 7.1
on page 309]{Bousfield-Kahn (1972)} which is an Atiyah-Hirzebruch type spectral
sequence and comes from a skeletal filtration. We will
need some preliminaries for its construction.

For every non-negative integer $p$, define the category $[p]$ whose objects are
$\{0,1,2,\ldots , p\}$,  with precisely one morphism from $i$ to $j$ if $i \leq j$ and no
morphism otherwise.  
Let \boldDelta \ be the category of finite ordered sets, i.e. objects are
the categories $[0],[1],[2], \ldots$ and morphisms are the functors from $[p]$ to
$[q]$.   In other words, the morphisms from an object $\{0,1,2,\ldots , p\}$ to an object
$\{0,1,2,\ldots , q\}$ are the monotone increasing functions.  A {\em simplicial set} is a
contravariant
$\boldDelta$-set.  There is a covariant $\boldDelta$-space $\Delta_{\bullet}$ which sends an
object $[p]$ to the standard $p$-simplex.  The {\em geometric realization of a
simplicial set $X_{\bullet}$} is the space 
$|X_{\bullet}| = X_{\bullet} \otimes_{\boldDelta} \Delta_{\bullet}$.
Recall that the {\em nerve of a category
$\calc$} is the simplicial set
$$N_p\calc ~ = ~ \func([p],\calc)$$ and its classifying space $B\calc$ is the
geometric realization $|N_{\bullet}\calc|$ of its nerve.
Next we introduce a similar construction. Denote by
$\widetilde{\func}([p],\calc)$ the equivalence classes of covariant
functors from $[p]$ to $\calc$, where two such functors are called equivalent if they are
related by a natural transformation whose evaluation at
any object is an isomorphism. More explicitly,
$\func([p],\calc)$ consists of the set of diagrams in $\calc$ of the
shape
$$c_0 ~ \xrightarrow{\phi_0} ~ c_1 ~ \xrightarrow{\phi_1} ~ c_2 ~
\xrightarrow{\phi_2} ~  \cdots ~ \xrightarrow{\phi_{p-1}} ~ c_p$$ and
$\widetilde{\func}([p],\calc)$ is the
set of equivalence classes under the following equivalence relation on
$\func([p],\calc)$: Two such diagrams
$(c_{\ast},\phi_{\ast})$ and
$(c_{\ast}^{\prime},\phi_{\ast}^{\prime})$ are equivalent if there
is a commutative diagram with
isomorphisms as vertical maps
\begin{center} $
\begin{CD} c_0 & @>\phi_0>> & c_1 &  @>\phi_1 >> &
\cdots &  @>\phi_{p-1} >> & c_p\\ @V \cong VV & & @V \cong VV & & & &
& & @V \cong VV\\ c_0^{\prime} &
@>\phi_0^{\prime} >> & c_1^{\prime} &  @>\phi_1^{\prime} >> & \cdots &
@>\phi_{p-1}^{\prime} >>& c_p^{\prime}
\end{CD}$\end{center} Define $\widetilde{N}_{\bullet}\calc$ to be the
simplicial set given by
$$\widetilde{N}_p\calc ~ = ~ \widetilde{\func}([p],\calc)$$ and
$\widetilde{B}\calc$ to be its geometric
realization.  We shall proceed to develop basic properties of $\widetilde B \calc$
analogous to those of $B \calc$, as discussed in \cite[pp. 227-229]{Davis-Lueck
(1998)}.

If $\calc$ is a groupoid, i.e. all morphisms are isomorphisms, then
$\widetilde{B}\calc$ is the (discrete) set of
isomorphism classes of objects. If $\calc$ is a category such that the
identity morphisms are the only
isomorphisms in $\calc$, then $B\calc$ is the same as $\widetilde{B}\calc$.
For instance, $\widetilde{B}[1] = B[1] = [0,1]$.

\begin{lemma} \label{lem: tilde classifying space and products} The
projections from $\calc_0 \times \calc_1$ to $\calc_i$
for $i=0,1$ induce a homeomorphism
$$\widetilde{B}(\calc_0 \times \calc_1) ~
\to \widetilde{B}\calc_0 \times \widetilde{B}\calc_1.$$
\end{lemma}
\proofl Given two simplicial sets $A_{\bullet}$ and $B_{\bullet}$, define their product
$A_{\bullet} \times B_{\bullet}$ by sending  $[p]$ to the product $A_p \times B_p$.
The projections induce a homeomorphism
\cite[page 43]{Lamotke (1968)}
$$|A_{\bullet} \times B_{\bullet}| ~ \to |A_{\bullet}| \times |B_{\bullet}|.$$ Now the claim follows
since the projections induce
isomorphisms of simplicial sets
$$\widetilde{\func}([p],\calc_0 \times \calc_1) ~
\to ~
\widetilde{\func}([p],\calc_0) \times
\widetilde{\func}([p],\calc_1). \qedl$$

Given two objects $?$ and $??$ in $\calc$, define the category
$\cata{?}{\calc }{??}$ as follows: An object is a diagram
$$
? \xrightarrow{\alpha} c \xrightarrow{\beta} ??
$$
 in $\calc$. A
morphism from
$? \xrightarrow{\alpha} c \xrightarrow{\beta} ??$ to
$? \xrightarrow{\alpha^{\prime}} c^{\prime} \xrightarrow{\beta^{\prime}} ??$ 
is a commutative diagram in $\calc$ of
the shape
\begin{center} $\begin{CD} ?  & @>\alpha >> & c & @>\beta >> & ??\\
@V\id VV & & @V\phi VV & & @V\id VV \\ ?  &
@>\alpha^{\prime} >> & c^{\prime} & @>\beta^{\prime} >> & ??
\end{CD}$\end{center} Let $\widetilde{B}_p(\cata{?}{\calc }{??})$
be the $p$-skeleton of $\widetilde{B}(\cata{?}{\calc }{??})$. We will regard
$\widetilde{B}(\cata{?}{\calc }{??})$ as a contravariant
$\calc \times \calc^{\op}$-space where $?$ is the variable in $\calc$ and
$??$ the variable in $\calc^{\op}$. Since
maps on classifying spaces induced by functors are cellular, we get a
filtration of the contravariant $\calc \times
\calc^{\op}$-space
$\widetilde{B}(\cata{?}{\calc }{??})$ by the contravariant $\calc
\times \calc^{\op}$-spaces
$\widetilde{B}_p(\cata{?}{\calc }{??})$ such that
$$\widetilde{B}(\cata{?}{\calc }{??}) ~ = ~
\colim_{p \to \infty} \widetilde{B}_p(\cata{?}{\calc }{??}).$$ Let
$\overline{\mor(?,??)}$ be the category whose
set of objects is
$\mor(?,??)$ and whose only morphisms are the identity morphisms of
objects. Consider the functor
$$\pr \colon  \cata{?}{\calc }{??} ~ \to \overline{\mor(?,??)}
\hspace{10mm} \left( ? \xrightarrow{\alpha} c \xrightarrow{\beta}
??\right)
 ~ \mapsto \left(\beta \circ \alpha \colon  ?
 \to ?? \right).$$ It induces a map of contravariant $\calc \times \calc^{\op}$-spaces
$$\widetilde{B}(\pr) \colon  \widetilde{B}(\cata{?}{\calc }{??}) ~
\to \widetilde{B}(\overline{\mor(?,??)}) = \mor(?,??).$$ Let
$X$ be a contravariant $\calc$-space. We
obtain contravariant
$\calc$-spaces
$X \otimes_{\calc} \widetilde{B}(\cata{?}{\calc }{??})$ and
$X \otimes_{\calc} \mor(?,??)$ where the
tensor product is taken over the variable $??$. Define a map of
contravariant
$\calc$-spaces
\begin{eqnarray}
& p \colon  X \otimes_{\calc} \widetilde{B}(\cata{?}{\calc }{??}) 
~ \xrightarrow{\id \otimes_{\calc} \widetilde{B}(\pr)} ~
X \otimes_{\calc} \mor(?,??)  ~ \xrightarrow{\cong} ~ X, &
\label{weak homotopy equivalence p}
\end{eqnarray}
where the second map is the canonical isomorphism given by
$x \otimes \phi \mapsto X(\phi)(x)$.

Let $n.d.\widetilde{N}_p(\cata{?}{\calc }{??})$ denote the set of
non-degenerate $p$-simplices of the
simplicial set $\widetilde{N}_{\bullet}(\cata{?}{\calc}{??})$. Elements are
given by classes of diagrams
$$? \xrightarrow{\alpha} c_0 ~ \xrightarrow{\phi_0} ~ c_1 ~
\xrightarrow{\phi_1} ~ c_2 ~ \xrightarrow{\phi_2} ~
\cdots ~ \xrightarrow{\phi_{p-1}} ~ c_p \xrightarrow{\beta} ??$$ such
that no $\phi_i$ is an isomorphism.

\begin{lemma} \label{lem: basic properties of tilde B cata}
\begin{enumerate}

\item \label{lem: basic properties of tilde B cata: colimits}
We have
$$X \otimes_{\calc} \widetilde{B}(\cata{?}{\calc }{??}) ~ = ~
\colim_{p \to \infty} X \otimes_{\calc}\widetilde{B}_p(\cata{?}{\calc
}{??})$$ as contravariant $\calc$-spaces;

\item  \label{lem: basic properties of tilde B cata: pushouts}
There is a pushout of contravariant $\calc$-spaces whose vertical
maps are $(p-1)$-connected cofibrations of
contravariant $\calc$-spaces
\begin{center}
$\begin{CD}
\left(X \otimes_{\calc} n.d.\widetilde{N}_p(\cata{?}{\calc }{??}) \right)
\times S^{p-1} & \to & X \otimes_{\calc}
\widetilde{B}_{p-1}(\cata{?}{\calc }{??})
\\ @V\id \otimes_{\calc} \inc VV   @V\id \otimes_{\calc} \inc VV
\\
\left(X \otimes_{\calc} n.d.\widetilde{N}_p(\cata{?}{\calc }{??}) \right)
\times D^p & \to & X \otimes_{\calc}
\widetilde{B}_p(\cata{?}{\calc }{??};
\end{CD}$
\end{center}

\item \label{lem: basic properties of tilde B cata: weak homotopy equivalence}
The map
$p \colon X \otimes_{\calc} \widetilde{B}(\cata{?}{\calc }{??}) \to  X$ 
defined in \eqref{weak homotopy equivalence p}
is a weak homotopy equivalence of contravariant
$\calc$-spaces.
\end{enumerate}
\end{lemma}
\proofl 
\eqref{lem: basic properties of tilde B cata: colimits} 
Since the functor $X \otimes_{\calc} -$ from the
category of covariant $\calc$-spaces to the
category of spaces has a right adjoint, it is compatible with
colimits.\\[2mm] 
\eqref{lem: basic properties of tilde B cata: pushouts} There is a canonical
$CW$-structure on the geometric realization of a simplicial set whose
cells are in bijective correspondence with
the non-degenerate simplices \cite[page 39]{Lamotke (1968)}. Hence we
get  the following pushout of contravariant
$\calc \times \calc^{\op}$-spaces
\begin{center}
$\begin{CD} n.d.\widetilde{N}_p(\cata{?}{\calc }{??}) \times S^{p-1} &
\to
&\widetilde{B}_{p-1}(\cata{?}{\calc }{??})
\\ @VVV  @VVV
\\ n.d.\widetilde{N}_p(\cata{?}{\calc }{??}) \times D^p & \to
&\widetilde{B}_p(\cata{?}{\calc }{??}).
\end{CD}$
\end{center} Since the functor $X \otimes_{\calc } -$ has a right
adjoint, it is compatible with pushouts.
For any space $Z$ there is a natural homeomorphism
$$\left(X \otimes_{\calc} n.d.\widetilde{N}_p(\cata{?}{\calc }{??}) \right)
\times Z \to X \otimes_{\calc}
\left(n.d.\widetilde{N}_p(\cata{?}{\calc }{??})
\times Z\right)$$ This shows that the diagram appearing in assertion 
\eqref{lem: basic properties of tilde B cata: pushouts} is
a pushout of contravariant
$\calc$-spaces. As the inclusion of $S^{p-1}$ into $D^p$ is a $(p-1)$-connected
cofibration of spaces, the left
vertical and hence also the right vertical arrow are $(p-1)$-connected
cofibrations of contravariant
$\calc$-spaces.
\\[2mm] 
\eqref{lem: basic properties of tilde B cata: weak homotopy equivalence}
Fix an object $c$ in $\calc$. Define a functor
$$j^{\prime} \colon  \overline{\mor(c,??)} \to
\cata{c}{\calc }{??} \hspace{10mm}
 \left(c \xrightarrow{\alpha} ??\right) ~\mapsto ~ \left(c
 \xrightarrow{\id} c \xrightarrow{\alpha} ?? \right).$$ It induces a map of
spaces by
$$j\colon  X(c) \xrightarrow{\cong} X \otimes_{\calc}
\widetilde{B}(\overline{\mor(c,??)})
\xrightarrow{\id_X \otimes_{\calc} \widetilde{B}j^{\prime}} X
\otimes_{\calc} \widetilde{B}(\cata{c}{\calc }{??}).$$
Define a natural transformation
$$S \colon  j' \circ \pr(c,??) \to \id_{c\da {\calc}\da {??}}$$ by assigning to an object
$c \xrightarrow{\alpha} d \xrightarrow{\beta}  ??$ in
$\cata{c}{\calc }{??}$ the morphism in $\cata{c}{\calc }{??}$
\begin{center}
$\begin{CD} c  @>\id >> c @> \beta \circ \alpha >> ?? \\ @V\id VV 
@V\alpha VV   @V\id VV \\ c @> \alpha >> d
@>\beta >> ??
\end{CD}$
\end{center} Thus we have  a homotopy of maps of covariant $\calc$-spaces
$$h^{\prime} \colon  \widetilde{B}(\cata{c}{\calc }{??}) \times [0,1]
\to
\widetilde{B}(\cata{c}{\calc }{??}) \times \widetilde{B}[1]
\to
\widetilde{B}(\cata{c}{\calc }{??} \times [1]) \to
\widetilde{B}(\cata{c}{\calc }{??})$$ where the first map comes from the
identification
$[0,1] = \widetilde{B}[1]$, the second from the homeomorphism
of Lemma \ref{lem: tilde classifying space and
products} and the third from interpreting $S$ as a functor
$$S \colon  \cata{c}{\calc }{??} \times [1] \to
\cata{c}{\calc }{??}.$$ It induces a homotopy of maps of spaces
$$h = \id_X \otimes_{\calc} h^{\prime}\colon  X \otimes_{\calc}
\widetilde{B}(\cata{c}{\calc }{??}) \times [0,1]
\to X \otimes_{\calc} \widetilde{B}(\cata{c}{\calc }{??}).$$
One easily checks that
$p(c) \circ j$ is the identity on $X(c)$ and $h$ is a homotopy from
$j \circ p(c)$ to the identity on
$X \otimes_{\calc} \widetilde{B}(\cata{c}{\calc }{??})$. Hence
$$p(c) \colon  X \otimes_{\calc} \widetilde{B}(\cata{c}{\calc }{??})
\to  X(c)$$ is a homotopy equivalence and in particular a weak
homotopy equivalence for all objects
$c$. This finishes the proof of Lemma \ref{lem: basic properties of tilde B
cata}. \qedl

\begin{remark} \em
\label{rem: p is  not a homotopy equivalence} Notice that the map
$p \colon  X \otimes_{\calc} \widetilde{B}(\cata{?}{\calc }{??}) \to  X$
 is {\em not} a homotopy equivalence of contravariant
$\calc$-spaces. In the proof of Lemma
\ref{lem: basic properties of tilde B cata} we have constructed a homotopy
inverse and a corresponding homotopy for
$p(c)$ for each object $c$, but they do not fit together to an homotopy
inverse of $p$ as a map of contravariant
$\calc$-spaces. Therefore it is important that we use the (co-)homology
$H(X;\bfE)$ which satisfies the WHE-axiom. \em \end{remark}

We can now apply \cite[Theorem 4.7]{Davis-Lueck (1998)} to the filtration of
$X \otimes_{\calc} \widetilde{B}(\cata{?}{\calc }{??})$
 by the subspaces
$X \otimes_{\calc} \widetilde{B}_p(\cata{?}{\calc }{??})$  and obtain:

\begin{theorem} \label{the: new spectral  sequence} Let $X$ be a contravariant
$\calc$-space and $\bfE$ a covariant
respectively contravariant $\calc$-spectrum. Let $H_p^{\calc}(X;\bfE)$
respectively
$H^p_{\calc}(X;\bfE)$ be the associated homology respectively cohomology
theories satisfying the WHE-axiom.

\begin{enumerate}

\item There is a spectral (homology) sequence
$(E^r_{p,q},d^r_{p,q})$ whose $E^1$-term is given by
$$E^1_{p,q} ~ = ~ H_q^{\calc}(X \otimes_{\calc} n.d.\widetilde{N}_p
(\cata{?}{\calc }{??});\bfE).$$ The first
differential
$$d^1_{p,q}\colon  ~ H_q^{\calc}(X \otimes_{\calc}
n.d.\widetilde{N}_p(\cata{?}{\calc }{??});
\bfE) ~ \to ~ H_q^{\calc}(X \otimes_{\calc}
n.d.\widetilde{N}_{p-1} (\cata{?}{\calc}{??});\bfE)
$$ is
$\sum_{i=0}^p  (-1)^i \cdot H_q^{\calc}(\id \otimes_{\calc} d^p_i;
\bfE)$ where
$d_i^p \colon  n.d.\widetilde{N}_p(\cata{?}{\calc }{??})
\to n.d.\widetilde{N}_{p-1} (\cata{?}{\calc }{??})$ is the $i$-th face map.
This spectral sequence converges to
$H^{\calc}_{p+q}(X;\bfE)$, i.e. there is an ascending filtration
$F_{p,m-p}H^{\calc}_{m}(X,\bfE)$ of $H^{\calc}_{m}(X,\bfE)$ such that
$$F_{p,q}H^{\calc}_{p+q}(X,{ \bf E}) /F_{p-1,q+1}H_{p+q}^{\calc}(X,{ \bf E})
~ \cong ~ E^{\infty}_{p,q};$$

\item Assume that $\bfE$ is a $\calc$-$\Omega$-spectrum, i.e. for each
object $c$ the spectrum $\bfE(c)$ is an
$\Omega$-spectrum. Then there is a spectral (cohomology) sequence
$(E_r^{p,q},d_r^{p,q})$ whose $E^1$-term is given by
$$E_1^{p,q} ~ = ~ H^q_{\calc}(X \otimes_{\calc} n.d.\widetilde{N}_p(\cata{?}
{\calc }{??});\bfE).$$ The first
differential
$$d_1^{p,q}\colon  ~H^q_{\calc}(X \otimes_{\calc}
 n.d.\widetilde{N}_p(\cata{?}{\calc }{??});\bfE) ~ \to ~
 H_{\calc}^q(X \otimes_{\calc} n.d.\widetilde{N}_{p+1}
(\cata{?}{\calc}{??});\bfE)
$$ is
$\sum_{i=0}^{p+1} (-1)^i \cdot H^q(\id \otimes_{\calc} d^p_i;\bfE)$. 
If one of the following conditions is satisfied:
\begin{enumerate}

\item The filtration is finite, i.e. there is an integer $n > 0$ such that for
any diagram
$$c_0 \xrightarrow{\phi_0} c_1 \xrightarrow{\phi_1}\cdots
\xrightarrow{\phi_{n-1}} c_n$$ one of the morphisms $\phi_i$ is an
isomorphism;

\item There is $n \in \zz$ such that $\pi_q(\bfE(c))$ vanishes for all
objects $c \in \Ob(\calc)$ and $q < n$;
\end{enumerate} 
then the spectral sequence converges to
 $H^{p+q}_{\calc}(X;\bfE)$, i.e. there is a  descending filtration
$F^{p,m-p}H^{m}_{\calc}(X,{ \bf E})$ of
$H^{m}_{\calc}(X,{ \bf E})$ such that
$$F^{p,q}H_{\calc}^{p+q}(X;\bfE)/ F^{p+1,q-1}H_{\calc}^{p+q}(X;{ \bf E}) ~
\cong ~ E_{\infty}^{p,q}.$$

\end{enumerate}

\end{theorem}

Next we want to analyse the $E^1$-term further. Let $\Is(\calc )$ be the
set of isomorphism classes $\overline{c}$
of objects $c$ in $\calc$. Fix for any isomorphism class $\overline{c}$ a
representative $c \in \overline{c}$. For two objects $c$ and $d$ let
$\mor_{\not\cong}(c,d)$ be the subset of $\mor(c,d)$ consisting
of all morphisms from $c$ to $d$ which are
not isomorphisms. For an element
$$\overline{c}_{\ast} = (\overline{c_0},\overline{c_1}, \ldots
,\overline{c_p})
\in \prod_{i=0}^p \Is(\calc )$$ and $p \ge 1$ define a
left-$\aut(c_p)$-right-$\aut(c_{0})$-set
\begin{eqnarray}
S(\overline{c}_{\ast}) & = & \mor_{\not\cong}(c_{p-1},c_p)
\times_{\aut (c_{p-1})}
\mor_{\not\cong}(c_{p-2},c_{p-1}) \times_{\aut (c_{p-2})} \nonumber
\\ & & \hspace{40mm} \ldots \times_{\aut(c_{1})} \mor_{\not\cong}(c_{0},c_1).
\label{definition of S(overline(c)_*)}
\end{eqnarray}
If $A$ is a right-$\aut(c_p)$-set and $B$
is a left-$\aut(c_0)$-set, define
$A \times_{\aut (c_p)} S(\overline{c}_{\ast}) \times_{\aut (c_0)} B$ 
in the obvious way for $p \ge 1$ and by
$A \times_{\aut(c_0)} B$ for $p=0$. One easily checks that the map
$$\coprod_{\substack{\overline{c}_{\ast} \in \prod_{i=0}^p \Is( \calc ),\\
S(\overline{c}_*) \not= \emptyset}}~
\mor(c_p,??) \times_{\aut (c_p)} S(\overline{c}_{\ast})
\times _{\aut (c_0)} \mor(?,c_0) ~ \to ~
n.d.\widetilde{N}_p(\cata{?}{\calc }{??})$$ which sends the
element represented by
$$(\phi_p,\phi_{p-1}, \ldots, \phi_0,\phi)
\in \mor(c_p,??) \times \mor_{\not\cong}(c_{p-1},c_p)
\times \ldots \times \mor_{\not\cong}(c_0,c_1)
\times \mor(?,c_0)$$ to the class of
$$? \xrightarrow{\phi} c_0 \xrightarrow{\phi_0} c_1
\xrightarrow{\phi_1} \cdots
\xrightarrow{\phi_{p-1}} c_p \xrightarrow{\phi_p} ??$$ is natural in $?$
and $??$ and bijective. Since there is a
natural isomorphism of contravariant $\calc$-spaces
\begin{multline*}
X \otimes_{\calc} \mor(c_p,??) \times_{\aut (c_p)}
S(\overline{c}_{\ast})
\times _{\aut (c_0)} \mor(?,c_0) 
\\ \xrightarrow{\cong} ~  X(c_p)
\times_{\aut(c_p)} S(\overline{c}_{\ast})
\times _{\aut (c_0)} \mor(?,c_0)
\end{multline*}
 we conclude:

\begin{lemma} \label{lem: computations of E_1 terms} There are
identifications for the $E^1$-terms of the spectral
sequences in Theorem \ref{the: new spectral  sequence}
$$E^1_{p,q} ~ = ~ \bigoplus_{\substack{\overline{c}_{\ast}
\in \prod_{i=0}^p \Is( \calc ),\\
S(\overline{c}_*) \not= \emptyset}}~ H_q^{\calc}(X(c_p) \times_{\aut(c_p)}
S(\overline{c}_{\ast})
\times_{\aut(c_0)} \mor(?,c_0);\bfE)$$ and
$$E_1^{p,q} ~ = ~ \bigoplus_{\substack{\overline{c}_{\ast} \in
\prod_{i=0}^p \Is( \calc ),\\S(\overline{c}_*) \not= \emptyset}}~
 H^q_{\calc}(X(c_p) \times_{\aut(c_p)}
S(\overline{c}_{\ast})
\times_{\aut(c_0)} \mor(?,c_0);\bfE). $$
\end{lemma}

If the category satisfies an additional condition, we can do
a much better job identifying the $E^1$-terms.

\begin{definition} \label{def: free category} We call $\calc$ {\em left-free}
 if for any two objects $c$ and $c'$ the left $\aut(c')$-action on
 $\mor(c,c')$ given by composition is free. 
\end{definition}

For any group $G$ and any family $\calf$ of subgroups the orbit category
$\Or(G,\calf)$ is left-free since any $G$-map of homogeneous $G$-spaces
$G/H \to G/K$ is surjective. 

Let $Z$ be a (left) $G$-space and let $\bfF$ be spectrum
with an action of $G$ by maps of spectra.
We can interpret $\bfF$ also as a covariant $\Or(G,\caltr)$-spectrum, where $\caltr$
is the family consisting of the trivial subgroup of $G$. We write
\begin{eqnarray}
H_q^G(Z;\bfF) & := & H_q^{\Or(G,\caltr)}(Z,\bfF).
\label{H_q^G(Z;bfF)  :=  H_q^{Or(G,caltr)}(Z,bfF)}
\end{eqnarray}
Explicitly we get after a choice of a free $G$-$CW$-complex $Z'$ together with a $G$-map
$u\colon Z' \to Z$ which is a weak equivalence after forgetting the group action
\begin{eqnarray}
H_q^G(Z;\bfF) & = & \pi_q(Z'_+\wedge_G \bfF).
\label{H_q^G(Z;bfF)  =  pi_q(Z'wedge_G bfF)}
\end{eqnarray}
For instance $H_q^G(\point;\bfF) = \pi_q(EG_+ \wedge_G \bfF)$.

\begin{lemma} \label{lem: computations of E_1 terms for free categories}
Suppose $\calc$ is left-free. Then there are
identifications for the
$E^1$-terms of the spectral sequences in Theorem \ref{the: new spectral sequence}
$$E^1_{p,q} ~ = ~ \bigoplus_{\substack{\overline{c}_{\ast}
\in \prod_{i=0}^p \Is( \calc ),\\S(\overline{c}_*) \not= \emptyset}}~ 
H_q^{\aut(c_0)}(X(c_p) \times_{\aut(c_p)}
S(\overline{c}_{\ast}) ;\bfE(c_0))$$ and
$$E_1^{p,q} ~ = ~ \bigoplus_{\substack{\overline{c}_{\ast}
\in \prod_{i=0}^p \Is( \calc ),\\S(\overline{c}_*) \not= \emptyset}}~ 
H^q_{\aut(c_0)}(X(c_p) \times_{\aut(c_p)}
S(\overline{c}_{\ast});\bfE(c_0))$$
where $X(c_p) \times_{\aut(c_p)} S(\overline{c}_{\ast})$ means
$X(c_0)$ for $p =0$. 

With this identification the first differential can be written as
$$d^1_{p,q} = \sum_{i=0}^p(-1)^i(d^1_{p,q})_i$$
for certain maps
\begin{multline*}
(d^1_{p,q})_0 \colon  H_q^{\aut(c_0)}(X(c_p) \times_{\aut(c_p)}S(\overline{c}_{\ast}) ;\bfE(c_0)) \\
\to  
 H_q^{\aut(c_1)}(X(c_p) \times_{\aut(c_p)}S(\overline{c}_1, \ldots, \overline{c}_{p}) ;\bfE(c_1)),
\end{multline*}
\begin{multline*}
(d^1_{p,q})_i \colon  H_q^{\aut(c_0)}(X(c_p) \times_{\aut(c_p)}S(\overline{c}_{\ast}) ;\bfE(c_0))\\
\to   
H_q^{\aut(c_0)}(X(c_p) \times_{\aut(c_p)}S(\overline{c}_0, \ldots  ,\overline{c}_{i-1},
\overline{c}_{i+1} , \ldots ,
\overline{c}_{p}) ;\bfE(c_0)),
\end{multline*}
for $0< i <p$, and
\begin{multline*}
(d^1_{p,q})_p \colon  H_q^{\aut(c_0)}(X(c_p) \times_{\aut(c_p)}S(\overline{c}_{\ast}) ;\bfE(c_0))\\
\to   
 H_q^{\aut(c_0)}(X(c_{p-1}) \times_{\aut(c_{p-1})}S(\overline{c}_0, \ldots, \overline{c}_{p-1})
;\bfE(c_0)).
\end{multline*}
For $0 < i \leq p$, these maps are induced by maps of $\aut(c_0)$-sets given by concatenation. 
The description of $(d^1_{p,q})_0$ is more difficult, due to the change of group and coefficients. 
Let  $c_0 \to c_1$ denote the full subcategory of $\calc$ with objects
$\{c_0,c_1\}$.  We label the inclusions of categories 
$i \colon\aut(c_0) \to (c_0 \to c_1)$,  $j \colon \aut(c_1) \to (c_0 \to c_1)$, 
and $k \colon  (c_0 \to c_1) \to \calc$.  Let $Y$ be the $\aut(c_1)$-space
$$
X(c_p) \times_{\aut(c_p)} S(\overline{c}_1, \ldots, \overline{c}_{p}).
$$
Then with the above identification of $E^1$ the map $(d^1_{p,q})_0$ is given by
\begin{multline*}
H^{\aut(c_0)}_q(Y \times_{\aut(c_1)} \mor_{\not\cong}(c_0,c_1); i^* k^* \bfE) 
\\ 
\to  H^{c_0\to c_1}_q(i_*(Y \times_{\aut(c_1)} \mor_{\not\cong}(c_0,c_1));  k^* \bfE)
\\
\to H^{c_0\to c_1}_q(j_*Y;  k^* \bfE) 
\xrightarrow{\cong} H^{\aut(c_1)}_q(Y;j^*k^*\bfE),
\end{multline*}
where the first map is the map $\Phi_i$ from Lemma 
\ref{lem: induction} \eqref{lem: induction: map Phi_F},  
the middle map is given
by the map of $(c_0 \to c_1)$-spaces given as the adjoint of the inclusion of $\aut(c_0)$-spaces
$$
Y \times_{\aut(c_1)} \mor_{\not\cong}(c_0,c_1) \to Y \times_{\aut(c_1)} \mor(c_0,c_1) = i^*j_*Y,
$$
and the last map is $(\Phi_j)^{-1}$.

Similar statements are valid in cohomology.
\end{lemma}

\proofl For an object $c_0$  of $\calc$, let 
$i_{c_0} \colon  \aut(c_0) \to \calc$ be the corresponding
inclusion of categories.  Then, according to 
Lemma \ref{lem: computations of E_1 terms}, the $E^1$-term of the 
$p$-chain spectral sequence is identified with 
$$
E^1_{p,q} ~ = ~ \bigoplus_{\substack{\overline{c}_{\ast}
\in \prod_{i=0}^p \Is( \calc ),\\S(\overline{c}_*) \not= \emptyset}}~
 H_q^{\calc}((i_{c_0})_*(X(c_p) \times_{\aut(c_p)}
S(\overline{c}_{\ast}));\bfE)
$$
The maps $\Phi_{i_{c_0}}$ send the sum below to the sum above
to
$$ \bigoplus_{\substack{\overline{c}_{\ast}
\in \prod_{i=0}^p \Is( \calc ),\\S(\overline{c}_*) \not= \emptyset}}~ 
H_q^{\aut(c_0)}(X(c_p) \times_{\aut(c_p)}
S(\overline{c}_{\ast}) ;\bfE(c_0))$$
This map between sums is an isomorphism  by Lemma  \ref{lem: induction}
\eqref{lem: induction: map Phi_F}, since $\calc$ is left-free.
We have thus established the first identification in the lemma.  
The  identifications of the differentials are established similarly.
\qedl

\begin{remark} \em
\label{rem: reduction in spectral sequence} Lemma
\ref{lem: computations of E_1 terms for free categories}
shows for a left-free category $\calc$ what the spectral sequence does.
Namely, it reduces the computation of the
(co-)homology groups of spaces and spectra over a
category to the special case of spaces and
spectra over a group.   The most important case of the $p$-chain spectral sequence is where
$X=\point$ is the constant functor given by a point at every object.  In this case the
$\aut(c_0)$-sets are all discrete and hence a disjoint union of homogeneous spaces $\aut(c_0)/H$,
and so the differentials
$(d^1_{p,q})_i$ for $0 < i \leq p$ all involve change of group maps $BH_0 \to BH_1$.   The remaining
differential
$(d^1_{p,q})_0$ is more subtle and should be
thought of as some sort of assembly map.  We will comment further on this map in the next
section.\em
\end{remark}

\begin{remark} \label{rem: EI-categories} \em We have put some
effort into avoiding  the assumption that $\calc$
is an {\em EI-category}, i.e that all endomorphisms are isomorphisms.
Otherwise we would have excluded the orbit
category
$\Or(G,\calf)$ for $\calf$ the family of virtually cyclic subgroups of $G$
(see \cite[Example 1.32]{Lueck (1989)}).
But this category appears in the Isomorphism Conjecture in algebraic
$K$ or $L$-theory of Farrell-Jones. The Baum-Connes Conjecture,
however, uses only the orbit category
$\Or(G,\calfin)$ for $\calfin$ the family of finite subgroups of $G$ which
is an EI-category. 

If one has an EI-category, the bookkeeping simplifies a little bit. The EI-property makes
it possible to define a partial ordering on $\Is(\calc)$ by
$$\overline{c} \le \overline{d} ~ \Longleftrightarrow
\mor(c,d) \not= \emptyset.$$ We write $\overline{c} < \overline{d}$ if
$\overline{c} \le \overline{d}$ and $\overline{c}
\not= \overline{d}$ holds. A $p$-chain $\overline{c_{\ast}}$ in $\calc$ is a
sequence
$$\overline{c_0} < \overline{c_1} < \ldots < \overline{c_p}.$$ Let
$ch_p(\calc)$ be the set of $p$-chains. Now one
can replace the index set
$$\left\{\left.\overline{c}_* \in \prod_{i=0}^p \Is(\calc)\right|~ 
 S(\overline{c}_*) \not= \emptyset\right\}$$
 in Lemma \ref{lem: computations of E_1 terms} and
Lemma \ref{lem: computations of E_1 terms for free categories} by
$ch_p(\calc)$ and replace in the definition of
$S(\overline{c_{\ast}})$ the set $\mor_{\not\cong}(c_{i-1},c_i)$ by
$\mor(c_{i-1},c_i)$.
\em  \end{remark}

\begin{example} \label{exa: orbit category}\em Let 
$\calf$ be a family of subgroups of the (discrete) group $G$. We get
a bijection
$$\left\{ (H) \mid H \in \calf \right\} ~
\xrightarrow{\cong} ~ \Is(\Or(G,\calf))
\hspace{10mm} (H) \mapsto \overline{G/H}$$ where $(H)$ denotes the
conjugacy class of a subgroup $H \subseteq G$.
Let $N\!H$ be the normalizer and $W\!H = N\!H/H$ be the Weyl group of $H
\subseteq G$. We obtain a bijection
$$W\!H ~ \xrightarrow{\cong} \aut(G/H)  \hspace{10mm} gH  \mapsto
\left(R_{g^{-1}} \colon  G/H \to G/H \right)$$ where $R_{g^{-1}}$
maps $g^{\prime}H$ to $g^{\prime}g^{-1}H$.
For
$$(\overline{G/H_i} \mid i =0,1, \ldots p) ~ \in ~ \prod_{i=0}^p \Is(\\Or(G,\calf))$$
we get
\begin{eqnarray*}
S(\overline{G/H_i} \mid i =0,1, \ldots p) 
&  = &
\map_{\not\cong}(G/H_{p-1},G/H_p)^G \times_{W\!H_{p-1}}\\
 & & \hspace{24mm} \ldots \times_{W\!H_1} \map_{\not\cong}(G/H_0,G/H_1)^G,
\end{eqnarray*}
where $\map_{\not\cong}(G/H_{i-1},G/H_i)^G$ is the set of
$G$-maps which are not bijective.

Suppose that $\Or(G,\calf)$ is an EI-category. Then a $p$-chain
$\overline{G/H_0} <  \ldots < \overline{G/H_p}$
is the same as a sequence of conjugacy classes of subgroups
$$(H_0) < \ldots  <  (H_p)$$ where $(H_{i-1}) < (H_i)$ means that
$H_{i-1}$ is subconjugated, but not
conjugated to $(H_i)$. The $W\!H_p$-$W\!H_0$-set associated to such a
$p$-chain is
\begin{eqnarray*}
S((H_0) < \ldots < (H_p)) & = & 
\map(G/H_{p-1},G/H_p)^G \times_{W\!H_{p-1}}
\\ & & \hspace{30mm} \ldots \times_{W\!H_1} \map(G/H_0,G/H_1)^G. 
\end{eqnarray*}
The $E^1$-terms of the
spectral sequences in Theorem \ref{the: new spectral  sequence} become
$$E^1_{p,q} ~ = ~ \bigoplus_{(H_0) < \ldots <  (H_p)}~
H_q^{W\!H_0}(X(G/H_p) \times_{W\!H_p} S((H_0) <\ldots <  (H_p)) ;\bfE(G/H_0))$$ 
and
$$E_1^{p,q} ~ = ~ \prod_{(H_0) < \ldots <  (H_p)}~
H^q_{W\!H_0}(X(G/H_p) \times_{W\!H_p} S((H_0) < \ldots  < (H_p)) ;\bfE(G/H_0))$$ 
where  $X(G/H_p) \times_{W\!H_p} S((H_0) < \ldots < (H_p))$
means $X(G/H_0)$ for $p =0$.
\em \end{example}

There is a module or chain complex version of the spectral sequence above.
In the sequel we use the notation and
language of \cite{Lueck (1989)}. Let $R$ be a commutative associative
ring with unit and $\calc$ be a small
category. One replaces the contravariant space $X$ by a contravariant
$R\calc$-chain complex $C$ which satisfies
$C_p = 0$ for $p < 0$ and the covariant respectively contravariant
spectrum
$\bfE$ by a covariant respectively contravariant $R\calc$-chain complex
$D$ which may have non-trivial chain
modules in negative dimensions. The role of contravariant $\calc$-$CW$-complexes is now
played by projective contravariant $R\calc$-chain complexes. The tensor product
$\otimes_{\calc}$ is replaced by the tensor product $\otimes_{R\calc}$
of
$R\calc$-modules and the mapping space is replaced by the $R$-module
of homomorphisms of $R\calc$-modules
$\hom_{R\calc}$. The homology
$H_p^{\calc}(X;\bfE)$ now becomes
$\Tor^{R\calc}_p(C,D)$ and the cohomology $H^p_{\calc}(X;\bfE)$ now
becomes
$\Ext_{R\calc}^p(C,D)$. Notice that a $R\calc$-module can be interpreted as
a chain complex concentrated in
dimension $0$. The proof of the next result is analogous to the results of
this section and generalizes the
spectral sequence in \cite[section 17]{Lueck (1989)}.

\begin{theorem} \label{the: module version} Let $M$ be a contravariant
$R\calc$-module and $N$ be a covariant
respectively contravariant $R\calc$-module. Suppose that $\calc$ is left-free. Then:

\begin{enumerate}

\item There is a spectral (homology) sequence
$(E^r_{p,q},d^r_{p,q})$ whose $E^1$-term is given by
$$E^1_{p,q} ~ = ~
\bigoplus_{\overline{c}_{\ast} \in
\prod_{i=0}^p \Is( \calc )}~
\Tor _q^{R[\aut(c_0)]}(M(c_p) \times_{R[\aut(c_p)]}
RS(\overline{c}_{\ast}) ,N(c_0))$$ where $R[\aut(c_i)]$ is the
group ring of $\aut(c_i)$ with coefficients in $R$ and
$RS(\overline{c}_{\ast})$ is the free $R$-module generated
by the set $S(\overline{c}_{\ast})$. The spectral sequence converges to
$\Tor^{R\calc}_{p+q}(M,N)$;

\item There is a spectral (cohomology) sequence $(E_r^{p,q},d_r^{p,q})$
whose $E^1$-term is given by
$$E_1^{p,q} ~ = ~ \prod_{\overline{c}_{\ast} \in
\prod_{i=0}^p \Is( \calc )}~
\Ext ^q_{R[\aut(c_0)]}(M(c_p) \times_{R[\aut(c_p)]}
RS(\overline{c}_{\ast});N(c_0)).$$

The spectral sequence converges to
$\Ext_{R\calc}^{p+q}(M,N)$. 
\end{enumerate}

\end{theorem}

\bigskip


\typeout{-------------------------- section 3  --------------------------}

\tit{Assembly Maps}{Assembly Maps}

Let $F \colon  \calb \to \calc$ and $\bfE \colon  \calc \to \SPECTRA$ be covariant functors.
We introduced their assembly map 
$$
H^\calb_q(\point;F^*\bfE) \xrightarrow{\Phi_F}    H^\calc_q(\point;\bfE)
$$
in Definition \ref{def: assembly map associated to a functor}.
Recall that we sometimes write $\bfE$ instead of $F^*\bfE$ to simplify notation.

Let $F\colon  \calb \to \calc$ be a covariant functor.  For any object 
$c \in \Ob C$, define the {\em undercategory}
$c \da  F$ and the {\em overcategory} $F \da  c$ as follows.  An object of  
$c \da  F$  is a pair $(b,\phi\colon c \to F(b))$   where
$b$ is an object in $\calb$ and $\phi$ a morphism in $\calc$.
A morphism $f$ from $(b,\phi \colon c \to F(b))$ to $(b',\phi' \colon c \to F(b'))$ 
is a morphism $f\colon b \to b'$ in $\calb$ satisfying $F(f) \circ \phi = \phi'$.
An object of  $F \da  c$  is a 
pair $(b,\phi \colon  F(b) \to c)$. A morphism
$f$ from $(b,\phi\colon F(b)\to c)$ to $(b',\phi'\colon F(b')\to c)$ 
is a morphism $f\colon b \to b'$ in $\calb$ satisfying $\phi' \circ F(f) = \phi$.
We denote the under and overcategories associated to the identity functor  $F : \calc \to \calc$ by $c
\da
\calc$ and
$\calc
\da c$.

A covariant functor $F\colon  \calb \to \calc$ is
{\em cofinal} if for every object $c$ of $\calc$, the classifying space $B(c\da F)$ is
contractible. 

For a category $\calc$, let $E\calc$ be the contravariant $\calc$-space given by
$$
E\calc(?) = B(? \da \calc).
$$
The  $\calc$-map  $E\calc \to \point$ is a $\calc$-$CW$-approximation by \cite[p. 230]{Davis-Lueck (1998)}.

\bigskip

\begin{theorem}[Cofinality Theorem] \label{the: Cofinality Theorem}  
Let $F \colon  \calb \to \calc$ be a cofinal covariant functor. Let
$\bf E \colon  \calc \to \SPECTRA $
be a covariant functor.  Then the assembly map
$$H^\calb_q(\point;F^*\bfE)\to  H^\calc_q(\point;\bfE)$$
is an isomorphism.  
\end{theorem}

\proofl The proof can be found in \cite[4.4]{Hollender-Vogt (1992)}. 
Here is a proof in our language. Notice that
$F_*E\calb(?)$ is a $\calc$-$CW$-complex by  Lemma \ref{lem:
induction} \eqref{lem: induction: CW-complexes}. 
Since there is a natural isomorphism
$F_*E\calb(?) \cong B(?\da  F)$ and $B(?\da  F)$ is contractible for
each $?$ by cofinality, the unique map $F_*E\calb(?) \to \point$ is a $\calc$-$CW$-approximation.  
Using Lemma \ref{lem: induction} \eqref{lem: induction: map Phi_F}, we conclude
$$H^\calb_q(\point; F^*\bfE) 
~ = ~
H^\calb_q(E\calb; F^*\bfE) 
~ = ~ 
H^\calc_q(F_*E\calb; \bfE)
~ = ~ 
H^\calc_q(\point; \bfE). \qedl
$$

For example, let $\calc$ be a category with a final object $c_0$. Let $\calb$ be the subcategory 
with single object $c_0$ and only the identity morphism. Then the inclusion 
functor $F\colon  \calb \to \calc$ is cofinal,  and so Theorem
\ref{the: Cofinality Theorem} shows that
$$\pi_q(\bfE(c_0))=H^{\calb}_q(\point;F^*\bfE) \cong H_q^\calc(\point;\bfE),$$ 
which is a well-known  fact about homotopy colimits.  This also follows from the observation
that $\point_{\calc} = \mor_{\calc}(?,c_0)$ is a $\calc$-$CW$-complex in this case.  Note that $\Or(G)$ has a final object $G/G$.

Given a functor $F : \calb \to \calc$ and an object $c$ in $\calc$, there is a commutatative square of 
functors
$$
\begin{CD} F \da c @>F_c>> \calc \da c \\
@VP_cVV @VQ_cVV  \\
\calb @>F>> \calc
\end{CD}
$$
where 
\begin{align*}
P_c(b, \phi : F(b) \to c) & = b, \\ 
F_c(b, \phi : F(b) \to c ) & =(F(b), \phi : F(b) \to c),\\ 
Q_c(c', \phi : c' \to c) &=
c'.
\end{align*}

Let $\bfE$ be a covariant $\calc$-spectrum.  Since $\calc\da  c$ has a final object, namely $(c,\id_c)$, we get an
identification
$$H_q^{\calc \da  c}(\point;Q_c^*\bfE) = \pi_q(\bfE(c)).$$

\begin{theorem}[Analogue of Quillen's Theorem A] \label{the: analogue of Quillen's Theorem A} 
Let $n$ be an integer and let $\calp$ be a set of primes.
Let $F \colon  \calb \to \calc$ and  $\bfE \colon  \calc \to \SPECTRA$ 
be covariant functors.
Suppose that  for all 
objects $c$ in $\calc$ the assembly map induced by $F_c$ 
$$H_q^{F \da c}(\point;F_c^*Q_c^*\bfE)\to H_q^{\calc \da  c}(\point;Q_c^*\bfE) = \pi_q(\bfE(c))$$ 
is a $\calp$-isomorphism for all $q \leq n$, where a $\calp$-isomorphism is a map of abelian
groups which becomes bijective after inverting all elements in $\calp$.

Then the assembly map induced by $F$
$$H^\calb_q(\point;F^*\bfE)  \to H^\calc_q(\point;\bfE)$$
is a $\calp$-isomorphism for all $q \leq n$.
\end{theorem}

\begin{remark} \label{QTA} \em
Quillen's Theorem A \cite{Quillen (1973)} states that for a covariant functor \break$F\colon \calb \to
\calc$, if for all objects $c$ of $\calc$, one has that 
$BF_c\colon  B(F \da c) \to B(C\da c) \simeq \point $ is a homotopy
equivalence, then $BF\colon  B\calb \to
B\calc$ is a homotopy equivalence.  
\em
\end{remark}

Theorems similar to and special cases of Theorem \ref{the: analogue of Quillen's Theorem A} are well-known.  Since we lack a
reference and for the reader's convenience,  we give a proof in our case.
We need some preparation. 

\begin{lemma} \label{lem: P(F downarrow c to calb)}
Let $X$ be a contravariant $(F\da  c)$-space. Then there is a natural homeomorphism of 
contravariant $\calb$-spaces
$$v\colon {P_c}_* X(?)  
\xrightarrow{\cong} 
\coprod_{\phi \in \mor_{\calc}(F(?),c)} ~ X(?,\phi\colon F(?) \to c).$$
In particular ${P_c}_*$ maps  weak equivalences of
contravariant $(F\da  c)$-spaces to  weak equivalences of
contravariant $\calb$-spaces.
\end{lemma}
\proofl A typical element in the domain of $v$ is represented by a pair
$(x,f)$, where $x \in X(b,\psi\colon F(b) \to c)$ and $f\colon ? \to b$ is a morphism
in $\calb$.  The image of it under $v$ is given by the image of $x$ under the map
$$X(b,\psi\colon F(b) \to c) \to X(?,\psi \circ F(f)\colon F(?) \to c)$$ 
which is given by $X$ applied to the morphism
$$f\colon  (?,\psi \circ F(f) \colon F(?) \to c) \to  (b,\psi\colon F(b) \to c)$$
in $F\da  c$. One easily checks that this is consistent with 
the tensor-relations appearing in the definition of
${P_c}_* X(?)$ as a tensor product. The map in the other direction sends
$x \in X(?,\phi\colon F(?) \to c)$ to the element in ${P_c}_*X$ represented by the pair
$(x,\id_{?})$. One easily checks using the tensor 
relation that this map is surjective and $v$ composed
with it is the identity on $\coprod_{\phi \in \mor_{\calc}(F(?),c)} ~ X(?,\phi\colon F(?) \to c)$.

Since a disjoint union of weak equivalences of spaces is again a weak equivalence of spaces,
$v$ sends weak equivalences to weak equivalences. \qedl

The contravariant $\calb$-space
${P_c}_*E(F\da  c)$ is actually a contravariant $\calc^{\op} \times \calb$-space,
where the covariant $\calc$-structure comes from the 
functoriality in the object $c$ in $\calc$. 

By a {\em homotopy $\calb$-$CW$-approximation of a $\calb$-space $Y$} we mean a weak homotopy equivalence
$Y' \to Y$ of $\calb$-spaces where $Y'$ has the homotopy type of a $\calb$-$CW$-complex. 

\begin{lemma} \label{homotopy CW} Let $F : \calb \to \calc$ be a covariant functor and let $X$ be a contravariant
$\calc$-$CW$-complex.  Then the map
$$
X \otimes_{\calc} {P_c}_*E(F\da  c) \to X \otimes_{\calc} {P_c}_*\point = X \otimes_{\calc} \mor_{\calc}(F(-),c)) = F^*X
$$
is a homotopy $\calb$-$CW$-approximation.
\end{lemma}

\proofl  Since $X$ is a $\calc$-$CW$-complex and, by Lemma \ref{lem: P(F downarrow c to calb)}, 
${P_c}_*E(F\da  c) \to  {P_c}_*\point$ is a weak homotopy equivalence, Theorem 3.11 of \cite{Davis-Lueck (1998)} shows
that $X \otimes_{\calc} {P_c}_*E(F\da  c) \to X \otimes_{\calc} {P_c}_*\point$ is weak homotopy equivalence.

To complete the proof we will show that for any contravariant $\calc$-$CW$-complex $Y$ and any covariant functor
$Z$ from $\calc$ into the category of contravariant 
$\calb$-$CW$-complexes that the contravariant $\calb$-space 
$Y \otimes_{\calc} Z$ has the homotopy type of a 
contravariant $\calb$-$CW$-complex. We first show by induction
that for each $n \ge -1$ the $\calb$-space $Y_n \otimes_{\calc} Z$ has the homotopy type of
a contravariant $\calb$-$CW$-complex. The induction beginning $n = -1$ is trivial because of 
$Y_{-1} \otimes_{\calc} Z = \emptyset$. In the induction step 
we can write $Y_n$ as a pushout of contravariant $\calc$-spaces
\comsquare{\coprod_{i \in I} S^n \times \mor_{\calc}(?,c_i)}{}{Y_n}
{}{}
{\coprod_{i \in I} D^{n+1} \times \mor_{\calc}(?,c_i)}{}{Y_{n+1}}
Since the functor $ - \otimes_{\calc} Z$ has a right adjoint, 
it is compatible with pushouts and we obtain a
pushout of contravariant $\calb$-spaces 
\comsquare{\coprod_{i \in I} S^n \times Z(c_i)}{}{Y_n \otimes_{\calc} Z}
{}{}
{\coprod_{i \in I} D^{n+1} \times Z(c_i)}{}{Y_{n+1} \otimes_{\calc} Z}
whose left vertical arrow is a cofibration of contravariant $\calb$-spaces. 
Since the contravariant spaces occuring in the left upper, right upper and left lower corner
have the homotopy type of contravariant $\calb$-$CW$-complexes, 
the same is true for $Y_{n+1} \otimes_{\calc} Z$.
Since each inclusion $Y_n \to Y_{n+1}$
is  a cofibration of contravariant $\calc$-$CW$-complexes, each inclusion
$Y_n \otimes_{\calc} Z \to Y_{n+1} \otimes_{\calc} Z$ is a cofibration of contravariant $\calb$-spaces 
which have the homotopy type of  contravariant $\calb$-$CW$-complexes.
Since $Y = \colim_{n \to \infty} Y_n$, then  $Y \otimes_{\calc} Z =  
\colim_{n \to \infty} Y_n \otimes_{\calc} Z$. Hence $Y \otimes_{\calc} Z$ and in particular
$X \otimes_{\calc} {P_c}_*E(F\da  c)$ have the homotopy type of a contravariant $\calb$-$CW$-complex.

 \qedl

Now we can prove Theorem 
\ref{the: analogue of Quillen's Theorem A}.\\
\proofl
We will factor the assembly map:
\begin{align*}
H^\calb_q(\point ; F^*\bfE) &\xleftarrow{\cong} H^\calb_q(F^*E\calc; F^*\bfE) & \text{w.h. invariance}\\
&\xleftarrow{\cong} \pi_q((E\calc \otimes_\calc {P_c}_*E(F\da c))\otimes_\calb F^*\bfE) & \text{Lemma \ref{homotopy CW}}\\
&= \pi_q(E\calc \otimes_\calc ({P_c}_*E(F\da c)\otimes_\calb F^*\bfE))   & \text{}\\
&= H^\calc_q(\point; {P_c}_*E(F\da c)\otimes_\calb F^*\bfE)   & \text{}\\
&= H^\calc_q(\point;E(F\da c)\otimes_{F\da c} {P_c}^*F^*\bfE) & \text{adjunction \eqref{adjunction (F_*,F^*) for
otimes}}\\ 
&= H^\calc_q(\point;E(F\da c)\otimes_{F\da c} {F_c}^*{Q_c}^*\bfE) & \text{}\\
&\to H^\calc_q(\point; \bfE) & \text{}
\end{align*}

By assumption the map 
$$
\pi_i(E(F\da c)\otimes_{F\da c} {F_c}^*{Q_c}^*\bfE) \to \pi_i(\bfE(c))
$$
is a $\calp$-isomorphism for each $i \le n$ and each object $c$ in $\calc$. 
Now Theorem  \ref{the: analogue of Quillen's Theorem A}
follows from a standard comparison argument applied to the spectral homology sequence
\cite[Theorem 4.7 (1)]{Davis-Lueck (1998)}.
 \qedl

Our motivation for studying assembly maps, as well as for creating the 
$p$-chain spectral sequence, is that the standard
conjectures of high-dimensional topology (Farrell-Jones and Baum-Connes) 
can be expressed in terms of assembly maps.  One
problem is that applying the $p$-chain spectral sequence directly to 
$\Or(G,\calvc)$  can be quite unpleasant, since
there are so many virtually cyclic subgroups.  We wish to make several remarks aimed 
at applying Theorems
\ref{the: Cofinality Theorem} and
\ref{the: analogue of Quillen's Theorem A} to the case where
$\calb =
\Or(G,\calf) \hookrightarrow \Or(G,\calg) = \calc$ where 
$\calf \subseteq \calg$ are families of subgroups of $G$.  The
aim is to ``simplify'' a family to make computations easier.  
In the sequel \emph{family} $\calf$ of subgroups of $G$ means a set
of subgroups of $G$ which is closed under conjugation.

An example is the following corollary which appeared $n = \infty$ and $\calp = \emptyset$ in
 \cite[Theorem A.10]{Farrell-Jones (1993a)} and for $\calp = \emptyset$  in
\cite[Theorem 2.3] {Lueck-Stamm (2000)}.

\begin{corollary}  
\label{cor: transitivity of assembly isomorphisms}
Let $\calf \subseteq \calg$ be families of subgroups of $G$.
Let $\bfE \colon  \Or(G, \calg) \to \SPECTRA$ be a covariant functor.  
For a subgroup $H$ of $G$, put $\calf_H = \{ K \mid  K \in \calf \text{ and } K \subseteq H\}$.  Let $n$ be an integer and 
let $\calp$ be a set of primes.
Suppose that for all $H \in \calg - \calf$
and for all $q \leq n$, the assembly map
$$
H^{\Or(H,  \calf_H)}_q(\point;\bfE) \to   H^{\Or(H,\calall)}_q(\point;\bfE)
$$
is a $\calp$-isomorphism.  Then  the assembly map 
$$
H^{\Or(G,  \calf)}_q(\point;\bfE) \to   H^{\Or(G,\calg)}_q(\point;\bfE)
$$
is a $\calp$-isomorphism for all $q \leq n$.
\end{corollary}

\proofl  First note that the condition in the statement of the lemma holds for 
all $H \in \calg$, since for $H \in \calf$,
$\Or(H,\calf_H)$ has a final object $H/H$.

Then  Theorem \ref{the: analogue of Quillen's Theorem A} applies where 
$F : \calb = \Or(G,\calf) \hookrightarrow \calc = \Or(G,\calg)$, since for 
all $H \in \calg$, there is an equivalence of categories
\begin{align*}
\Or(H,\calf_H) &\to F\da G/H  \\
H/K &\mapsto   (G/K, gK \mapsto gH)  .   \qedl
\end{align*}

We next apply the Cofinality Theorem.

\begin{lemma}  \label{lem: cofinal business}    
Let $\calf \subseteq \calg$ be families of subgroups of $G$.  
Suppose that for every $H \in \calg - \calf$ there exists 
$K_H \in \calf$ with the properties that $H \subseteq K_H$ and that for any 
$L \in \calf$ with $H \subseteq L$ we have $L \subseteq K_H$,
in other words, $K_H$ is the  maximal element of the set $\{L \in \calf\mid H \subseteq L\}$.

Then the inclusion functor $I\colon \Or(G,\calf) \to \Or(G,\calg)$ is cofinal.
\end{lemma}
\proofl 
We have to show for $G/H \in \Or(G,\calg)$ that $B(G/H \da I)$
is contractible. Recall the well-known fact that $B\calc$ is contractible if
$\calc$ contains an initial or final object $c_0$. 
If $H$ belongs to $\calf$, the object $(G/H,\id\colon G/H\to G/H)$ is an initial object in
$G/H \da I$. It remains to show for $H \in \calg - \calf$ that there exists a final object
in $G/H\da I$.

Let $K_H$ be the  maximal element of the set $\{L \in \calf\mid H \subseteq L\}$.
Let $\pr_H \colon G/H \to G/K_H$ be the canonical projection. It suffices to show for
any $L \in \calf$ and $G$-map $\phi\colon G/H \to G/L$ that there is a $G$-map 
$\overline{\phi}\colon G/L \to G/K_H$ satisfying $\overline{\phi} \circ \phi = \pr_H$. 
Since $\phi$ is surjective, $\overline{\phi}$ is uniquely determined by this property
and one easily checks that the object $(G/H,\pr_H\colon G/H \to G/K_H)$ is final
in $G/H\da I$. Choose $\gamma \in G$ with $\phi(1H) = \gamma L$.
Then $H \subseteq \gamma L \gamma^{-1}$. By assumption 
$\gamma L \gamma^{-1} \subseteq K_H$. Hence
we can define a $G$-map
$\overline{\phi} \colon G/L \to G/K_H$ by sending $gL$ to $g\gamma^{-1}K_H$.
One easily checks $\overline{\phi} \circ \phi = \pr_H$.
\qedl

\begin{corollary}  \label{cor: unique element} 
Let $\calf \subseteq\calg$ be family of all subgroups of $G$.  Suppose that every element 
of $\calg - \calf$ is contained in
a unique element of $\calf$.  

Then the inclusion functor $I\colon \Or(G,\calf) \to \Or(G,\calg)$ is cofinal.
\end{corollary}

The next result presents a systematically way to replace the family $\calg$ by a smaller family $\calg'$.

\begin{corollary}  \label{cor: passage from calg to calg'}
Let $\calg$ be a family of subgroups of $G$. Let $\calg' \subset\calg$ be the
subfamily consisting of all maximal elements of $\calg$, together with all 
elements of $\calg$ which are contained in no or in more than one maximal elements of $\calg$.  

Then the inclusion functor $I\colon \Or(G,\calg') \to \Or(G,\calg)$ is cofinal.
\end{corollary}

Let $H \subseteq G$ be a normal subgroup of $G$. Let $\{H\}$ be the family consisting
of a single subgroup, namely $H$, and let $\calsub(H)$ be the family
of subgroups of $H$.  Both are indeed closed under conjugation, since $H$ is normal.
Let $\bfE$ be a covariant $\Or(G,\calsub(H))$-spectrum.
The inclusion functor $I\colon \Or(G,\{H\}) \to \Or(G,\calsub(H))$ is cofinal
by Corollary \ref{cor: unique element}. Hence the assembly map
$$H_q^{\Or(G,\{H\})}(\point;I^*\bfE) \xrightarrow{\cong}
H_q^{\Or(G,\calsub(H))}(\point;\bfE)$$
is bijective by Theorem \ref{the: Cofinality Theorem}. Let $Q$ be the quotient group $G/H$.
There is an obvious isomorphism $\Or(G,\{H\}) \to \Or(Q,\caltr)$. It induces an isomorphism
\begin{eqnarray*}
H_q^{\Or(G,\{H\})}(\point;I^*\bfE) & \xrightarrow{\cong} &
H_q^Q(\point;\bfE(G/H))
\end{eqnarray*}
where we equip $\bfE(G/H)$ with the obvious $Q$-action and
$H_q^Q(\point;\bfE(G/H))$ has been introduced in 
\eqref{H_q^G(Z;bfF)  :=  H_q^{Or(G,caltr)}(Z,bfF)}. Thus we rediscover
the isomorphism from \cite[Lemma 2.6]{Lueck-Stamm (2000)}
\begin{eqnarray}
H_q^{\Or(G,\calsub(H))}(\point;\bfE) & \xrightarrow{\cong} &
H_q^Q(\point;\bfE(G/H)).
\label{computation for calsub(H)}
\end{eqnarray}

In the sequel $1$ denotes the trivial group.

\begin{definition} \label{def: assembling along a normal subgroup}
Let $H$ be a normal subgroup of $G$ and let $Q$ be the quotient group.  Let 
$\bfE \colon \Or(G,\calsub(H)) \to \SPECTRA$ be a covariant functor.  
Then the assembly map associated to the inclusion
$\Or(G,\caltr) \to \Or(G,\calsub(H))$ gives using the identifications
\eqref{H_q^G(Z;bfF)  :=  H_q^{Or(G,caltr)}(Z,bfF)} and 
\eqref{computation for calsub(H)} a map
$$
H_q^G(\point;\bfE(G/1)) \to H_q^Q(\point;\bfE(G/H)).
$$
It is called the {\em partial assembly map}, assembling along $H$.
\end{definition}

This partial assembly map can be identified with the composition
\begin{multline*}
\pi_q\left(EG_+ \wedge_G \bfE(G/1)\right) \xrightarrow{\pi_q(\id \wedge_G \bfE(\pr))} 
\pi_q\left(EG_+ \wedge_G \bfE(G/H)\right) \\
 ~ =  \pi_q\left((H\backslash EG)_+ \wedge_Q \bfE(G/H)\right)
\xrightarrow{\pi_q(f_+ \wedge_Q \id)} \pi_q\left(EQ_+ \wedge_Q \bfE(G/H)\right),
\end{multline*}
where $\pr\colon G/1 \to G/H$ is the projection and $f\colon H\backslash EG \to EQ$ is the 
up to $Q$-homotopy unique $Q$-map.

We next wish to make an analysis of a few of the terms 
and differentials in the $p$-chain spectral sequence. 
Recall the Weyl group $W\!H = N\!H/H$. We leave the elementary proof of the next result
to the reader.

\begin{theorem}
\label{the: information about d^1}
Let $\bfE \colon  \Or(G) \to \SPECTRA$ be a covariant functor.  
Consider the $p$-chain spectral sequence when $X
= \point$.  Then 

\begin{enumerate}
\item \label{the: information about d^1: summand of E^1 of 0-chain}

The summand of $E^1_{0,q}$ corresponding to a 
$0$-chain of the form $\overline{G/H}$ is isomorphic to
$$H_q^{W\!H}(\point;\bfE(G/H)) = \pi_q(EW\!H_+ \wedge_{W\!H} \bfE(G/H)),$$
where the $W\!H$-action on  $\bfE(G/H)$ comes from the action of the automorphisms
of $G/H$ in $\Or(G)$;

\item \label{the: information about d^1: summand of E^1 of 1-chain}
The summand of $E^1_{1,q}$ corresponding to a $1$-chain 
of the form $(\overline{G/1},\overline{G/H})$ is isomorphic to
$$H_q^{N\!H}(\point;\bfE(G/1)) = \pi_q(EN\!H_+\wedge_{N\!H} \bfE(G/1)),$$
where the $N\!H$-action on $\bfE(G/1)$ is the restriction of the $G$-action which is given by
the identification $G = \aut_{\Or(G)}(G/1)$;

\item  \label{the: information about d^1: d^1_{1,q}}
The component of $d^1_{1,q}$ corresponding to omitting 
$\overline{G/1}$ from a $1$-chain of the form
$(\overline{G/1},\overline{G/H})$ can be identified with the  
partial assembly map (see Definition \ref{def: assembling along a normal subgroup}
)
$$H_q^{N\!H}(\point;\bfE(G/1))
\to H_q^{W\!H}(\point;\bfE(G/H)),$$
for $H \subset N\!H$ and $\bfE$ considered as covariant functor
$\Or(N\!H,\calsub(H)) \to \SPECTRA$ by sending $N\!H/L$ to $\bfE(G/L)$.
\end{enumerate}

\end{theorem}

If there are no $p$-chains for $p > 1$, then the $p$-chain spectral sequence
gives a long exact sequence.

\begin{corollary} \label{cor: long exact sequence if M(M) holds}
Let $\calf$ be a family of subgroups of $G$. Suppose that any $G$-map  $G/H \to G/K$ 
is bijective, provided $H,K \in \calf$, $H \not= 1$, $K \not= 1$.  Let $\Lambda$ be the 
set of conjugacy classes $(H)$ of non-trivial subgroups in $\calf$. Let
$\bfE\colon  \Or(G,\calf) \to \SPECTRA$ be a functor. Then there
is a long exact sequence
\begin{multline*}
\ldots \to \bigoplus_{(H)\in \Lambda} H^{N\!H}_q(\point; \bfE(G/1)) \\
\xrightarrow{pa,-i} 
\left(\bigoplus_{(H)\in \Lambda} H^{W\!H}_q(\point; \bfE(G/H))\right)  
\bigoplus H_q^G(\point; \bfE(G/1)) \\
\xrightarrow{h\oplus a} H_q^{\Or(G,\calf)}(\point; \bfE) \to  
\bigoplus_{(H)\in \Lambda} H^{N\!H}_{q-1}(\point; \bfE(G/1)) \xrightarrow{pa,-i} \ldots
\end{multline*}
Here the map $pa$ is a sum of partial assembly maps, $i$ is a change of group map
associated to the inclusion $N\!H \to G$,
and $a$ is the assembly map induced by the
inclusion of categories $I \colon \Or(G,\caltr) \to \Or(G,\calf)$ under the identification
$H_q^{\Or(G,\caltr)}(\point;I^*\bfE) =  H_q^G( \point; \bfE(G/1))$.
The $N\!H$-action on $\bfE(G/1)$ is the restriction of the $G = \aut_{\Or(G)}(G/1)$-action
and the $W\!H$-action on  $\bfE(G/H)$ comes from the identification
$W\!H = \aut_{\Or(G)}(G/H)$.
\end{corollary}

For the reader's convenience we restate the exact sequence appearing
Corollary \ref{cor: long exact sequence if M(M) holds} in more
familiar terms given by homotopy groups
\begin{multline*}
\ldots \to \bigoplus_{(H)} \pi_q(EN\!H_+ \wedge_{N\!H} \bfE(G/1)) \\
\xrightarrow{pa,-i} 
\left(\bigoplus_{(H)} \pi_q(EW\!H\wedge_{W\!H}\bfE(G/H))\right)  
\bigoplus \pi_q(EG_+ \wedge_G \bfE(G/1)) \\
\xrightarrow{h\oplus a} H_q^{\Or(G,\calf)}(\point; \bfE) \to  
\bigoplus_{(H)} \pi_{q-1}(EN\!H_+ \wedge_{N\!H} \bfE(G/1)) \xrightarrow{pa,-i} \ldots
\end{multline*}

An example of family $\calf$ appearing in Corollary \ref{cor: long exact sequence if M(M) holds}  
is the family $\calfin$ of finite subgroups of a group $G$, provided 
every non-trivial finite subgroup of $G$ has prime order. We will see other
examples in Section \ref{sec: Examples}.


\typeout{-------------------------- section 4  --------------------------}

\tit{Examples}{Examples}

In this section we illustrate how the methods of the preceeding sections can be used
to compute $K_q(\zz G)$, $L^{\langle -\infty \rangle}_q(\zz G)$ and $K_q(C^*_r(G))$ for
certain discrete groups $G$ explicitly.

Let $\calmfin$ be the subfamily of $\calfin$ consisting
of elements in $\calfin$ which are maximal in $\calfin$. 
Consider the following assertions concerning $G$:\\[1mm]
\begin{itemize}

\item[(M)] Every non-trivial finite subgroup of $G$ is contained in a unique maximal finite subgroup;

\item[(NM)] $M \in \calmfin ~ \Rightarrow ~ N\!M = M$;

\item[(VCL)] If $V$ is an infinite virtually cyclic subgroup of $G$, then
$V$ is isomorphic to $\zz$;

\item[(BC)] The Baum-Connes Conjecture for $G$ is true, i.e. the assembly map
$$H_q^{\Or(G,\calfin)}(\point;\bfK^{\topo}) \to 
H_q^{\Or(G)}(\point;\bfK^{\topo}) = K_q(C^*_r(G)) 
$$
is bijective for all $q \in \zz$;

\item[(FJK$_N$)]
The Isomorphism Conjecture of Farrell-Jones for algebraic $K$-theory
is true for $\zz G$ in the range $q \le N$ for a 
fixed element $N \in \zz \coprod \{\infty\}$,
i.e. the assembly map
$$
H_q^{\Or(G,\calvc)}(\point;\bfK^{\alg}) \to 
H_q^{\Or(G)}(\point;\bfK^{\alg}) = K_q(\zz G)
$$
is bijective for $q \in \zz$ with $q \le N$;

\item[(FJL)]
The Isomorphism Conjecture of Farrell-Jones for algebraic $L$-theory
is true for $\zz G$, i.e.
the assembly map
$$
H_q^{\Or(G,\calvc)}(\point;\bfL^{\langle -\infty \rangle}) 
\to H_q^{\Or(G)}(\point;\bfL^{\langle -\infty \rangle}) 
= L^{\langle - \infty \rangle}_q(\zz G)
$$
is bijective for $q \in \zz$;

\item[(FJL{[1/2]})]
The Isomorphism Conjecture of Farrell-Jones for algebraic $L$-theory
is true for $\zz G$ after inverting $2$ , i.e.
the assembly map
$$
H_q^{\Or(G,\calvc)}(\point;\bfL^{\langle -\infty \rangle}) 
\to H_q^{\Or(G)}(\point;\bfL^{\langle -\infty \rangle}) = L^{\langle -
  \infty \rangle}_q(\zz G)
$$
is bijective after inverting $2$ for all $q \in \zz$..

\end{itemize}

Denote by $\underline{E}G$ \emph{the classifying space for proper $G$-actions}. It is a
$G$-$CW$-complex such that all of its isotropy groups are finite and the fixed point set of
each finite subgroup of $G$ is contractible. It has the universal property
that for any $G$-$CW$-complex $X$ with all isotropy groups finite, there is a $G$-map from $X$ to $\underline{E}G$, and any
two $G$-maps from $X$ to $\underline{E}G$ are $G$-homotopic. This property characterizes $\underline{E}G$ up to
$G$-homotopy.

For a group $H$, define $\widetilde{K}_q(\zz H)$ to be the cokernel of the map
$K_q(\zz[1]) \to K_q(\zz H)$ induced by the inclusion of the trivial
subgroup $1$ into $H$. Define $\widetilde{L}^{\langle - \infty \rangle}(\zz H)$,
$\widetilde{K}_q(C_r^*(H))$ and $\widetilde{K}_q(BH)$ analogously.

\begin{theorem} \label{the: computations based on IC}
Let $\zz \subseteq A \subseteq \qq$ be a ring  such that the order of any 
finite subgroup of $G$ is invertible in $A$.  Let $\Lambda$ be the set
of conjugacy classes (H) of non-trivial subgroups belonging to $\calmfin$. Then: 
\begin{enumerate} 

\item \label{the: computations based on IC:  (M), (NM) and (BC)}
If $G$ satisfies (M), (NM) and (BC), then there are exact sequences of topological $K$-groups
\begin{multline*}
\ldots \to \bigoplus_{(H) \in \Lambda} K_q(BH) \to  
\left(\bigoplus_{(H) \in \Lambda} K_q(C^*_r(H)) \right) \bigoplus K_q(BG) 
\to 
\\
K_q(C^*_r(G)) \to \bigoplus_{(H) \in \Lambda} K_{q-1}(BH) \to \ldots 
\end{multline*}
and
$$
0 \to \bigoplus_{(H) \in \Lambda} \widetilde{K}_q(C^*_r(H))
 \to K_q(C_r^*(G)) \to K_q(G\backslash \underline{E}G) \to 0. $$
The maps $K_q(BH) \to K_q(C^*_r(H))$ and $K_q(BG) \to
K_q(C^*_r(G))$ are the assembly maps associated to $\caltr \subseteq \calall$.
The maps  $K_q(BH) \to K_q(BG)$, $K_q(C^*_r(H)) \to K_q(C_r^*(G))$
and $\widetilde{K}_q(C^*_r(H)) \to K_q(C_r^*(G))$ are induced by the inclusion $H \to G$.
The second sequence splits after applying $- \otimes_{\zz} A$;

\item \label{the: computations based on IC:  (M), (NM), (VCL) and (FJL)}
If $G$ satisfies (M), (NM), (VCL) and (FJL), then there are exact sequences 
of algebraic $L$-groups 
\begin{multline*}
\ldots \to \bigoplus_{(H) \in \Lambda} H_q(BH; \bfL(\zz)) 
\to  
\left(\bigoplus_{(H) \in \Lambda} L_q^{\langle -\infty\rangle}(\zz H) \right) 
\bigoplus H_q(BG ; \bfL(\zz)) 
\\
\to 
L_q^{\langle -\infty\rangle}(\zz G ) \to
\bigoplus_{(H) \in \Lambda} H_{q-1}(BH; \bfL(\zz)) \to  \ldots
\end{multline*}
and
\begin{multline*}
\ldots \to H_{q+1}(G\backslash \underline{E}G ; \bfL(\zz)) 
\to \bigoplus_{(H) \in \Lambda} \widetilde{L}_q^{\langle -\infty\rangle}(\zz H)
\\
\to L_q^{\langle -\infty\rangle}(\zz G) \to H_q(G\backslash \underline{E}G ; \bfL(\zz)) \to \ldots
\end{multline*}
where $\bfL(\zz) = \bfL^{-\langle \infty \rangle}(\zz)$ is the algebraic $L$-theory-spectrum 
associated to the ring $\zz$. 
The maps $H_q(BH ; \bfL(\zz)) \to L_q^{\langle - \infty \rangle}(\zz H)$ and 
$H_q(BG ; \bfL(\zz)) \to L_q^{\langle - \infty \rangle}(\zz G)$
are the assembly maps associated to $\caltr \subseteq \calall$. The maps
and $H_q(BH ; \bfL(\zz)) \to H_q(BG ; \bfL(\zz))$,
$L_q^{\langle -\infty\rangle}(\zz H) \to L_q^{\langle -\infty\rangle}(\zz G )$ and
$\widetilde{L}_q^{\langle -\infty\rangle}(\zz H) \to L_q^{\langle -\infty\rangle}(\zz G )$
are  induced by the inclusion $H \to G$.
The second sequence splits after applying $- \otimes_{\zz} A$, more precisely
$$L_q^{\langle -\infty\rangle}(\zz G) \otimes_{\zz} A 
\to H_q(G\backslash \underline{E}G ; \bfL(\zz)) \otimes_{\zz} A$$
is a split-surjective map of $A$-modules;

\item \label{the: computations based on IC:  (M), (NM) and [(FJL[1/2])}
If $G$ satisfies (M), (NM), and (FJL[1/2]), then the conclusion of
assertion \eqref{the: computations based on IC:  (M), (NM), (VCL) and (FJL)} still holds
if we invert $2$ everywhere. Moreover,
the second sequence reduces to a short exact sequence
$$
0 \to \bigoplus_{(H) \in \Lambda} \widetilde{L}_q(\zz H)\left[\frac{1}{2}\right]
 \to L_q(\zz G)\left[\frac{1}{2}\right]
 \to H_q(G\backslash \underline{E}G ; \bfL(\zz)\left[\frac{1}{2}\right]
 \to 0. 
$$
which splits after applying 
$- \otimes_{\zz\left[\frac{1}{2}\right]}A\left[\frac{1}{2}\right]$.

(Recall that the decorations in $L$-theory do not matter after inverting $2$);

\item \label{the: computations based on IC:  (M), (NM), and (FJK_N)}
If $G$ satisfies (M), (NM),  and (FJK$_N$), then there is for
$q \in \zz, q \le N$ an isomorphism of Whitehead groups
$$\bigoplus_{(H) \in \Lambda} \Wh_q(H) \xrightarrow{\cong} \Wh_q(G),$$
where $\Wh_q(H) \to \Wh_q(G)$ is induced by the inclusion $H \to G$.

\end{enumerate}

\end{theorem}
\proofl
\eqref{the: computations based on IC:  (M), (NM) and (BC)}
The inclusion functor $\Or(G,\caltr \cup \calmfin) \to \Or(G,\calfin)$ is cofinal
because of assumption (M) and Corollary \ref{cor: passage from calg to calg'}.
Lemma \ref{the: Cofinality Theorem} implies  that
the assembly map associated to the inclusion $I\colon \Or(G,\caltr \cup \calmfin) \to \Or(G,\calfin)$ 
$$H_q^{\Or(G,\caltr \cup \calmfin)}(\point;\bfK^{\topo}) ~
\xrightarrow{\cong} H_q^{\Or(G,\calfin)}(\point;\bfK^{\topo})$$
is bijective for $q \in \zz$. Now the first long exact sequence
follows from assumptions (NM) and (BC)  and Corollary \ref{cor: long exact sequence if M(M) holds}.

From \cite[Lemma 2.8 (c)]{Lueck-Stamm (2000)} we obtain an exact sequence
\begin{multline*}
\ldots \to K_{q+1}(G\backslash \underline{E}G) \to 
H_q^{\Or(G,\calfin)}(\point;\widetilde{\bfK}^{\topo}) 
\\ 
\to H_q^{\Or(G,\calfin)}(\point;\bfK^{\topo}) \to  
K_q(G\backslash \underline{E}G) \to \ldots
\end{multline*}
such that the map
$$H_q^{\Or(G,\calfin)}(\point;\bfK^{\topo}) \otimes_{\zz} A ~ \to  ~
K_q(G\backslash \underline{E}G) \otimes_{\zz} A$$
is a split surjective maps of $A$-modules.
Here $\widetilde{\bfK}^{\topo}$ is a covariant $\Or(G)$-spectrum
satisfying $\pi_q(\widetilde{\bfK}^{\topo}(G/H)) = \widetilde{K}_q(C_r^*(H))$.
Since $\pi_q(\widetilde{\bfK}^{\topo}(G/1))$ vanishes for all $q \in \zz$,
we obtain from Corollary \ref{cor: passage from calg to calg'}, 
Corollary \ref {cor: long exact sequence if M(M) holds} and
assumption (NM) an isomorphism
$$\bigoplus_{(H) \in \Lambda} \widetilde{K}_q(C^*_r(H)) \xrightarrow{\cong}
H_q^{\Or(G,\calfin)}(\point;\widetilde{\bfK}^{\topo})$$
Thus we obtain a long exact sequence which splits after applying 
 $- \otimes_{\zz} A$
$$
\ldots \to K_{q+1}(G\backslash \underline{E}G) 
 \to \bigoplus_{(H) \in \Lambda} \widetilde{K}_q(C^*_r(H))
\to K_q(C_r^*(G)) \to K_q(G\backslash \underline{E}G) \to \ldots 
$$
Since  $\widetilde{K}_q(C^*_r(H))$ is a finitely generated torsionfree abelian group,
assertion \eqref{the: computations based on IC:  (M), (NM) and (BC)} follows.
\\[2mm]
\eqref{the: computations based on IC:  (M), (NM), (VCL) and (FJL)}
The proof is analogous to the one of assertion
\eqref{the: computations based on IC:  (M), (NM) and (BC)} except that
one additionally has to prove that the assembly map associated to the inclusion
$\Or(G,\calfin) \to \Or(G,\calvc)$ induces an isomorphism
$$H_q^{\Or(G,\calfin)}(\point;\bfL^{\langle - \infty \rangle}) 
\xrightarrow{\cong} 
H_q^{\Or(G,\calfin)}(\point;\bfL^{\langle - \infty \rangle}).$$
Because of Corollary \ref{cor: transitivity of assembly isomorphisms}
and assumption (VCL) it suffices to check that the following
assembly map is an isomorphism
$$H_q^{\Or(\zz,\caltr)}(\point;\bfL^{\langle - \infty \rangle}) 
\xrightarrow{\cong} 
H_q^{\Or(\zz)}(\point;\bfL^{\langle - \infty \rangle})$$ 
This follows from the Shaneson splitting  \cite[Theorem
5.1]{Shaneson(1969)}, the Rothenberg sequence and the fact that
$\widetilde{K}_i(\zz[\zz])$ for $i \le 0$ and $\Wh(\zz)$ vanish.
\\[2mm]
\eqref{the: computations based on IC:  (M), (NM) and [(FJL[1/2])}
The proof is analogous to the one of assertion 
\eqref{the: computations based on IC: (M), (NM), (VCL) and (FJL)} using the conclusion from
\cite{Yamasaki (1987)} that for any virtually cyclic group $V$ the assembly map
$$H_q^{\Or(V,\calfin)}(\point;\bfL^{\langle - \infty  \rangle})
\left[\frac{1}{2}\right] 
 \xrightarrow{\cong} 
H_q^{\Or(V)}(\point;\bfL^{\langle - \infty \rangle})
\left[\frac{1}{2}\right]$$ 
is bijective, and the conclusion from \cite[Proposition 22.34 on page 253]{Ranicki (1992)}
that $L_q(\zz H)\left[\frac{1}{2}\right]$ is a torsionfree
$\zz\left[\frac{1}{2}\right]$- module for a finite group $H$.
\\[2mm]
\eqref{the: computations based on IC:  (M), (NM), and (FJK_N)}
There are canonical identifications \cite[Lemma 2.4]{Lueck-Stamm (2000)}
\begin{eqnarray*}
H_q^{\Or(G)}(\point_{\Or(G)},\point_{\caltr};\bfK^{\alg})
& = &
\left\{\begin{array}{lcl}
\Wh_q(G) & & q \ge 2;
\\
\Wh(G) = \Wh_1(G) & & q = 1;
\\
\widetilde{K}_0(\zz G) = \Wh_0(G) & & q = 0;
\\
K_q(\zz G) = \Wh_q(G) & & q \le -1,
\end{array}
\right.
\end{eqnarray*}
where $\point_{\caltr} \subseteq \point_{\Or(G)}$ is the contravariant sub-$\Or(G)$-space
which sends $G/1$ to the one-point space and $G/H$ for $H \not= 1$ to the
empty set. The assembly map 
$$H_p^{\Or(V,\calfin)}(\point;\bfK^{\alg})
\xrightarrow{\cong} 
H_p^{\Or(V)}(\point;\bfK^{\alg})$$
is an isomorphism for all $q \in \zz$ if $V$ is $\zz$ or $\zz/2 \times \zz/2$
(see \cite[Lemma 2.5]{Lueck-Stamm (2000)}).
Let $V \subseteq G$ be an infinite virtually cyclic subgroup.
Any infinite virtually cyclic group $V$ admits an epimorphism with finite kernel $F$ to
$\zz$ or $\zz/2 \ast \zz/2$. Let $M$ be a maximal finite subgroup containing $F$.
For $g \in N\!F$ we have $F \subseteq M \cap gMg^{-1}$. If $F$ is non-trivial,
assumption  (M) and (NM) imply $g \in N\!M = M$ which contradicts $V \subseteq N\!F$. This shows
that $V$ is isomorphic to $\zz$ or $\zz/2 \ast \zz/2$. 
Hence the assembly map
$$H_q^{\Or(G,\caltr \cup \calmfin)}(\point;\bfK^{\alg}) \xrightarrow{\cong}
H_q^{\Or(G)}(\point;\bfK^{\alg}) = K_q(\zz G)$$
is an isomorphism for $q \le N$ by assumption (FJK$_N$) and
Corollary \ref{cor: transitivity of assembly isomorphisms}. This implies by a Five-Lemma
argument that there is for $q \le N$ an isomorphism
$$H_q^{\Or(G,\caltr \cup \calmfin)}(\point_{\Or(G,\caltr \cup \calmfin)},\point_{\caltr};\bfK^{\alg})
\xrightarrow{\cong} \Wh_q(G).$$
Now apply the $p$-chain spectral sequence (for pairs) to the source of
this map. Notice that here the $E^1$-term consists only of contributions by $0$-chains
associated to elements $H \in \calmfin$.
This finishes the proof of Theorem 
\ref{the: computations based on IC}. \qedl

\begin{remark} \label{rem: other coefficients} 
\em
The claims  for the topological $K$-theory of the complex reduced $C^*$-algebra
carry directly over to the real reduced $C^*$-algebra.

One can also replace in the computations of algebraic $K$ and $L$-theory the ring
$\zz$ by some commutative ring $R$ with unit. If $R$ contains $\qq$, then the assembly map
$$H_p^{\Or(V,\calfin)}(\point;\bfK^{\alg}) \to 
H_p^{\Or(V)}(\point;\bfK^{\alg}) = K_q(RV)$$
is bijective for any virtually cyclic group $V$, essentially because
$RV$ is a regular ring for any virtually cyclic group $V$ and $\qq \subseteq R$  and hence all Nil-terms
vanish. By Corollary \ref{cor: transitivity of assembly isomorphisms} the assembly map
$$H_p^{\Or(G,\calfin)}(\point;\bfK^{\alg}) \to 
H_p^{\Or(G,\calvc)}(\point;\bfK^{\alg})$$
is bijective and  one is for the computation of $K_q(RG)$ in the same fortunate 
situation as in $K_p(C^*_r(G))$, where it suffices to consider $\calfin$ instead of $\calvc$.
\em
\end{remark}

Now we give some groups for which Theorem \ref{the: computations based on IC} applies:

\begin{itemize}

\item Extensions $1 \to \zz^n \to G \to F \to 1$ for finite $F$ such that the conjugation
  action of $F$ on $\zz^n$ is free outside $0 \in \zz^n$. \\[1mm]
The conditions (M), (NM),  (BC) and (FJK$_1$) and (FJL[1/2])  are satisfied by
\cite[Theorem 1.2, Lemma 6.1 and Lemma 6.3]{Lueck-Stamm (2000)}. Hence the conclusions
appearing in assertions 
\eqref{the: computations based on IC:  (M), (NM) and (BC)},
\eqref{the: computations based on IC:  (M), (NM) and [(FJL[1/2])} and
\eqref{the: computations based on IC:  (M), (NM), and (FJK_N)} for $N = 1$ 
of Theorem \ref{the: computations based on IC} holds for $G$.
This has already been proved in \cite[Theorem 0.2]{Lueck-Stamm (2000)}, where further
information is given. The $L$-groups of $F$ are computed in \cite[Remark 6.4]{Lueck-Stamm (2000)}
for all decorations.  There a term UNIL appears which has meanwhile been
computed by Connolly and Davis \cite{Connolly-Davis (2002)};

\item Fuchsian groups $F$ \\[1mm]
The conditions (M), (NM),  (BC) and (FJK$_1$) and (FJL[1/2])  are satisfied for $F$
(see for instance \cite[Theorem 1.2 and Lemma 4.5]{Lueck-Stamm (2000)}). Hence the conclusions
appearing in assertions 
\eqref{the: computations based on IC:  (M), (NM) and (BC)},
\eqref{the: computations based on IC:  (M), (NM) and [(FJL[1/2])} and
\eqref{the: computations based on IC:  (M), (NM), and (FJK_N)} for $N = 1$ hold for $F$.
The $L$-groups of $F$ are computed in \cite[Remark 4.10]{Lueck-Stamm (2000)}
for all decorations.  There a term UNIL appears which has meanwhile been
computed by Connolly and Davis  \cite{Connolly-Davis (2002)}.   The computation of
$\Wh_q(F)$ for a Fuchsian groups $F$ for $q \le 1$ was independently carried out 
in \cite{Berkhove-Juan-Pineda-Pearson(2001)} and  \cite{Berkhove-Juan-Pineda-Pearson(2001b)} using the $p$-chain spectral
sequence. In  \cite{Lueck-Stamm (2000)} the larger class of cocompact planar
groups (sometimes also called cocompact NEC-groups) is treated.  The $p$-chain spectral sequence
was applied in \cite{Davis-Pearson (2002)} to prove the Gromov-Lawson-Rosenberg positive scalar 
curvature conjecture for manifolds whose fundamental group is Fuchsian;

\item One-relator groups $G$\\[1mm]
Let $G$ be a one-relator group. Let $ G = \langle (g_i)_{i \in I} \mid r \rangle$ 
be a presentation with one relation.  Then (BC) is satisfied
\cite{Beguin-Bettaib-Valette (1999)}.
We do not know whether (FJK$_N$) or (FJL) hold for $G$ in general, and we will assume in the
discussion below that they do. 

We begin with the case, where $G$ is torsionfree. Then there is a 
$2$-dimensional model for $BG$ which is given by the $2$-dimensional $CW$-complex
associated to any presentation with one generator
\cite[Chapter III \S\S 9 -11]{Lyndon-Schupp(1977)}. We obtain isomorphisms
\begin{eqnarray*}
K_p(BG) & \xrightarrow{\cong} & K_p(C^*_r(G));
\\
H_q(BG ; \bfK(\zz)) & \xrightarrow{\cong} & K_p(\zz G);
\\
H_q(BG ; \bfL(\zz)) & \xrightarrow{\cong} & L_p(\zz G);
\end{eqnarray*}

Now suppose that $G$ is not torsionfree.
Let $F$ be the free group with basis $\{g_i \mid i \in I\}$. Then $r$ is an element in
$F$. There exists an element $s \in F$ and an integer $m \ge 2$
such that $r = s^m$, the cyclic  subgroup $C$
generated by the class $\overline{s} \in G$ represented 
by $s$ has order $m$, any finite subgroup of $G$ is subconjugated to $C$ and 
for any $g \in G$ the implication
$g^{-1}Cg \cap C \not= 1 \Rightarrow g \in C$ holds.
These claims follows from
\cite[Propositions 5.17, 5.18 and 5.19 in II.5 on pages 107 and 108]{Lyndon-Schupp(1977)}.
Hence $G$ satisfies (M) and (NM) and up to conjugation there is precisely one maximal
finite subgroup, namely $C$.  Hence the
conclusion in assertion \eqref{the: computations based on IC:  (M), (NM) and (BC)} 
of Theorem \ref{the: computations based on IC} holds
for $G$ (see also \cite{Mislin (2002)}). 
If we assume (FJL[1/2]) or (FJK$_N$) respectively, this is also true
for assertions \eqref{the: computations based on IC:  (M), (NM) and [(FJL[1/2])},
and \eqref{the: computations based on IC:  (M), (NM), and (FJK_N)} of Theorem
\ref{the: computations based on IC}.

From now on suppose that (FJL) holds and that $m = |C|$ has odd order. Then
assertion \eqref{the: computations based on IC:  (M), (NM), (VCL) and (FJL)} of Theorem
\ref{the: computations based on IC} is true. Since the $\zz\left[\frac{1}{2}\right]$module
$\widetilde{L}_q^{-\langle  - \infty \rangle}(\zz C)\left[\frac{1}{2}\right]$
is finitely generated free, we get a short exact 
sequence
$$0 \to \widetilde{L}_q^{\langle  - \infty \rangle}(\zz C) \to 
L_q^{\langle  - \infty \rangle}(\zz G) \to 
H_q(G\backslash\underline{E}G ; \bfL(\zz)) \to 0,$$
which splits after inverting $m$. One easily checks that it induces 
a short exact  sequence
$$0 \to L_q^{\langle  - \infty \rangle}(\zz C) \to 
L_q^{\langle  - \infty \rangle}(\zz G) \to 
H_q(G\backslash\underline{E}G,\ast;\bfL(\zz)) \to 0,$$
which splits after inverting $m$.

We mention that 
there is a $2$-dimensional $CW$-model for $G\backslash\underline{E}G$ such that
there is precisely  one $0$-cell, precisely  one $2$-cell and there is a bijective correspondence
between the $1$-cells and the index set $I$. This follows from
\cite[Exercise 2 (c) II. 5 on page 44]{Brown (1982)}. The Atiyah-Hirzebruch spectral sequence
yields isomorphisms
$$H_q(G\backslash\underline{E}G,\ast;\bfL(\zz)) ~ \cong ~ 
\left\{
\begin{array}{lll}
H_2(G\backslash\underline{E}G;\zz/2) & & q = 0 ~ (4);
\\
H_1(G\backslash\underline{E}G;\zz)   & & q = 1 ~ (4);
\\
H_2(G\backslash\underline{E}G;\zz)   & & q = 2 ~ (4);
\\
H_1(G\backslash\underline{E}G;\zz/2) & & q = 3 ~ (4).
\end{array}
\right.$$
We mention without giving the proofs the following facts.
The abelian groups $H_2(G\backslash \underline{E}G;\zz)$ and 
$H_2(G\backslash \underline{E}G;\zz/2)$ vanish and there is an exact sequence
$0 \to H_1(C;\zz) \to H_1(BG;\zz) \to H_1(G\backslash \underline{E}G;\zz) \to 0$,
provided that the relation $r \in F$ does not belong to
the commutator $[F,F]$. We have 
$H_2(G\backslash \underline{E}G;\zz) \cong \zz $ and 
$H_2(G\backslash \underline{E}G;\zz/2)\cong \zz/2$, and 
$H_1(BG;\zz)) \cong H_1(G\backslash \underline{E}G;\zz)$ is a free abelian group of
rank $|I|-1$, provided that $r \in [F,F]$. 
\end{itemize}

There are a number of computations where $G$ is a two-dimensional crystallographic group.  The groups
$\Wh_q(G)$ for  $q \le 1$ were computed in \cite{Pearson (1998)}, the groups $K_0(\cc G)$ were computed in 
\cite{Pitkin (1995)}, and the groups
 $K_q(C^r_*(G))$ for $q \in \zz$ were
computed in \cite{Yang(1997)} (see also \cite{Pitkin (1995)} and \cite{Lueck-Stamm (2000)}).
These computations were done in a unified way in \cite[Section 5]{Lueck-Stamm (2000)}, as well as
 $\Wh_2(G)$ and 
$L_q(\zz G)\left[\frac{1}{2}\right]$ for $q \in \zz$.


\typeout{-------------------------- references  --------------------------}

Addresses\\[3mm] James F. Davis \hfill Wolfgang L\"uck \\ Department of
Mathematics \hfill Fachbereich f\"ur
Mathematik und Informatik\\ Indiana University \hfill Westf\"alische
Wilhelms-Universit\"at  \\ Blomington, IN
47405 \hfill 48149 M\"unster
\\ U.S.A. \hfill Einsteinstr. 62\\
\mbox{ }  \hfill Germany\\ email:  jfdavis@ucs.indiana.edu
\hfill lueck@topologie.math.uni-muenster.de\\ Fax-number: 812 855-0046
\hfill 0251 8338370
\\ http://www.indiana.edu/\verb+~+jfdavis/ \hfill
http://wwwmath.uni-muenster.de/
\\ \text{} \hfill math/u/lueck/

\end{document}